\documentclass[10pt]{amsart}
\usepackage{amssymb,amsmath,amsthm}
\usepackage{mathrsfs,dsfont,a4wide}
\theoremstyle{plain}
\newtheorem{theorem}{Theorem}[section]
\newtheorem{lemma}[theorem]{Lemma}
\newtheorem{proposition}[theorem]{Proposition}
\newtheorem{corollary}[theorem]{Corollary}
\newtheorem{remark}[theorem]{Remark}
\newtheorem{definition}[theorem]{Definition}
\theoremstyle{definition}
\theoremstyle{remark}
\numberwithin{equation}{section}

\newcommand{\as}{{\mathcal A}}
\newcommand{\hs}{{\mathcal H}}
\newcommand{\ks}{{\mathcal K}}

\newcommand{\gs}{{\mathcal G}}
\newcommand{\fs}{{\mathcal F}}

\newcommand{\ms}{{\mathcal M}}

\newcommand{\Es}{{\mathcal E}}

\newcommand{\Esup}{{\mathcal E}^s}
\newcommand{\Eb}{{\mathcal E}^b}

\newcommand{\R}{{\mathbb R}}
\newcommand{\N}{{\mathbb N}}

\newcommand{\Om}{\Omega}
\newcommand{\Omb}{\overline{\Omega}}

\newcommand{\weakst}{\stackrel{\ast}{\rightharpoonup}}

\newcommand{\weak}{\rightharpoonup}
\newcommand{\wlystar}{$\text{weakly}^*$\;}
\newcommand{\wstar}{$\text{weak}^*$\;}

\newcommand{\hn}{\hs^{N-1}}

\newcommand{\eps}{\varepsilon}
\newcommand{\tsub}{\,\tilde{\subseteq}\,}
\newcommand{\Sg}[2]{S^{#1}(#2)}

\newcommand{\res}{\mathop{\hbox{\vrule height 7pt width .5pt depth 0pt
\vrule height .5pt width 6pt depth 0pt}}\nolimits}

\title
[A $\Gamma$-convergence approach to stability of unilateral minimality properties]
{A $\Gamma$-convergence approach to stability \\
of unilateral minimality properties \\
in fracture mechanics and applications}
\author[A. Giacomini]
{Alessandro Giacomini}
\address[Alessandro Giacomini]{S.I.S.S.A., Via Beirut 2-4, 34014, Trieste,
Italy \& Max Planck Institute for
Mathematics in the Sciences, Inselstrasse 22, D-04103 Leipzig, Germany}
\email[A. Giacomini]{giacomin@sissa.it \& giacomin@mis.mpg.de}
\author[M. Ponsiglione]
{Marcello Ponsiglione}
\address[Marcello Ponsiglione]{S.I.S.S.A., Via Beirut 2-4, 34014, Trieste,
Italy \& Max Planck Institute for
Mathematics in the Sciences, Inselstrasse 22, D-04103 Leipzig, Germany}
\email[M. Ponsiglione (corresponding author)]{ponsigli@mis.mpg.de}
\begin{document}
\vskip .2truecm
\begin{abstract}
\small{
We prove the stability of a large class
of unilateral minimality properties which arise
naturally in the theory of crack propagation proposed by Francfort and
Marigo in \cite{FM}. 
Then we give an application to the quasistatic evolution of cracks
in composite materials.
}
\end{abstract}
\thanks{Corresponding author: Marcello Ponsiglione, Max Planck Institute for
Mathematics in the Sciences, Inselstrasse 22, D-04103 Leipzig, Germany, {\bf ponsigli@mis.mpg.de}}
\maketitle
%{\small \tableofcontents}

\section*{Introduction}
\label{intr}
In this paper we deal with the problem of stability of {\it unilateral minimality properties} with varying volume and surface energies, and we give an application to the study of crack propagation in composite materials.
\par
Let $K$ be a $(N-1)$-dimensional set contained in $\Om \subseteq \R^N$, and let $u$ be a possibly vector valued function on $\Om$ whose discontinuities are contained in $K$ and which is sufficiently regular outside $K$. We say that the pair $(u,K)$ is a {\it unilateral minimizer} with respect to the energy densities $f$ and $g$ if
\begin{equation}
\label{minuni}
\int_{\Om \setminus K}f(x,\nabla u(x))\,dx+\int_Kg(x,\nu)\,d\hn(x)
\le 
\int_{\Om \setminus H}f(x,\nabla v(x))\,dx+\int_Hg(x,\nu)\,d\hn(x).
\end{equation}
for every $(N-1)$-dimensional set $H$ containing $K$, and for every function $v$ 
whose discontinuities are contained in $H$ and which is sufficiently regular outside $H$. Here $\nu$ stands for the normal vector to $K$ and $H$ at the point $x$, while $\hn$ stands for the $(N-1)$-dimensional Hausdorff measure. $(u,K)$ is said to be {\it unilateral} minimizer because it is a minimum only among pairs $(v,H)$ with $H$ 
larger than $K$.
\par
The unilateral minimality property \eqref{minuni} is a key point in the theory of quasistatic crack evolution in elastic bodies proposed by Francfort and Marigo in \cite{FM} and which is inspired by the classical Griffith's criterion of crack propagation. In the framework of \cite{FM}, $\Om$ represents an hyperelastic body in the reference configuration, $u$ is its deformation, and $K$ represents a crack inside $\Om$ across which the deformation $u$ may jump. The total energy of the configuration $(u,K)$ is given by
\begin{equation}
\label{totener}
\Es(u,K):=\int_{\Om \setminus K}f(x,\nabla u(x))\,dx+\int_K g(x,\nu)\,d\hn(x).
\end{equation}
The first term is referred to as {\it bulk energy} of the body, while the second term is referred to as {\it surface energy} of the crack. The presence of $x$ in $f$ and $g$ takes into account possible inhomogeneities, while the presence of the normal $\nu$
in $g$ takes into account a possible anisotropy of the body.
\par
Following \cite{FM}, if $\Om$ is subject to a time dependent loading process, a quasistatic crack evolution can be described by a pair $(u(t),K(t))$ where the crack
$K(t)$ growths in time, $(u(t),K(t))$ satisfies the unilateral minimality property \eqref{minuni} at each time $t$, and the total energy \eqref{totener} evolves in relation with the power of external loads in such a way that no dissipation occurs.
\par
The unilateral minimality property \eqref{minuni} can be interpreted as a static equilibrium property along the irreversible process of crack growth. In fact an immediate consequence of \eqref{minuni} is that $u(t)$ is the elastic deformation in $\Om \setminus K(t)$ associated to the external load. As for the crack $K(t)$, \eqref{minuni} states a minimality condition only among enlarged cracks (unilateral minimality), taking thus into account the irreversibility of the process. Together with non dissipation, and under some regularity assumptions on the cracks, the unilateral minimality property implies that the Griffith's criterion is satisfied along the evolution (see \cite{DMT}).
\par
In \cite{FM} Francfort and Marigo suggest that the quasistatic evolution $(u(t),K(t))$ during the loading process can be obtained as a limit of a discretized in time evolution $(u_n(t),K_n(t))$ which by construction satisfies at each time the unilateral minimality property \eqref{minuni}. We are thus led to a problem of {\it stability} for unilateral minimizers, i.e. if the minimality property \eqref{minuni} is conserved in the passage from $(u_n(t),K_n(t))$ to $(u(t),K(t))$.
\par
The first mathematical result of stability for unilateral minimality properties was obtained by Dal Maso and Toader \cite{DMT} in a two dimensional setting under a topological restriction on the admissible cracks. They consider compact cracks with a bound on the number of connected components, and converging with respect to the Hausdorff metric.
An extension of this result for unilateral minimality properties involving the symmetrized 
gradients of planar elasticity has been done by Chambolle in \cite{Ch}, while an extension to higher order minimality properties in connection to quasistatic crack growth in a plate has been proved by Acanfora and Ponsiglione in \cite{AcP}.
\par
A second result of stability for unilateral minimality properties was obtained by Francfort and Larsen in \cite{FL}, where they give an existence result for quasistatic crack evolutions in the context of $SBV$ functions. In the framework of {\it generalized antiplanar shear} (i.e. $\Om \subseteq \R^N$, $N \ge 2$), the authors consider cracks $K$ which are rectifiable sets in $\Omb$, and associated displacements $u$ in $SBV(\Om)$ with jump set $S(u)$ contained in $K$. 
A key point for their result is the stability for unilateral minimizers of the form $(u_n,S(u_n))$ with bulk energy given by $f(x,\xi)=|\xi|^2$ and surface energy given by $g(x,\nu) \equiv 1$. More precisely, writing the minimality property in the equivalent form
\begin{equation*}
\label{minuni2}
\int_\Om |\nabla u_n|^2\,dx \le
\int_\Om |\nabla v|^2\,dx+\hn(S(v) \setminus S(u_n))
\qquad
\text{for all }v \in SBV(\Om)
\end{equation*}
(which corresponds to \eqref{minuni} with $H=S(u_n) \cup S(v)$), they prove that if $u_n \weak u$ weakly in $SBV(\Om)$ (see Section \ref{prel} for a definition), then $u$ satisfies the same minimality property. The main tool for proving stability is 
a geometrical construction which they called Transfer of Jump Sets \cite[Theorem 2.1]{FL}.
\par
The case in which $S(u_n)$ is replaced by a rectifiable set $K_n$ has been treated by Dal Maso, Francfort and Toader in \cite{DMFT}, where they consider also a Carath\'eodory bulk energy $f(x,\xi)$  quasiconvex and with $p$ growth assumptions in $\xi$, and a Borel surface energy $g(x,\nu)$ bounded and bounded away from zero.
They employ a variational notion of convergence for rectifiable sets which they called $\sigma^p$-convergence to recover a crack $K$ in the limit (see Section \ref{sigmaconv}), and they prove a Transfer of Jump Sets theorem for $(K_n)_{n \in \N}$ satisfying $\hn(K_n) \le C$ \cite[Theorem 5.1]{DMFT} in order to prove that minimality is preserved.
\par
In this paper we provide a different approach to the problem of stability of unilateral minimizer based on $\Gamma$-convergence which will permit also to treat the case of varying bulk and surface energy densities $f_n$ and $g_n$. We restrict our analysis to the scalar case. Our approach is based on the observation that the problem has a variational character. In fact, considering for a while the case of fixed energy densities $f$ and $g$ with $f$ convex in $\xi$, we have that if $(u_n,K_n)$ is a unilateral minimizer for the energy \eqref{totener},
then $u_n$ is a minimum for the functional
\begin{equation*}
\label{intrfuncn}
\Es_n(v):=\int_\Om f(x,\nabla v(x))\,dx+\int_{S(v) \setminus K_n}g(x,\nu)\,d\hn(x).
\end{equation*}
Then the problem of stability of unilateral minimizers can be treated in the framework of $\Gamma$-convergence which ensures the convergence of minimizers.
In Section \ref{gammasec}, using an abstract representation result by 
Bouchitt\'e, Fonseca, Leoni and Mascarenhas \cite{BFLM}, we prove that the $\Gamma$-limit (up to a subsequence) of the functional $\Es_n$ can be represented as
\begin{equation*}
\label{intrfunc}
\Es(v):=\int_\Om f(x,\nabla v(x))\,dx+\int_{S(v)}g^-(x,\nu)\,d\hn(x),
\end{equation*}
where $g^-$ is a suitable function defined on $\Om \times S^{N-1}$ determined only by $g$ and $(K_n)_{n \in \N}$, and such that $g^- \le g$. If we assume that $u_n \weak u$ weakly in $SBV(\Om)$, then by $\Gamma$-convergence we get that $u$ is a minimizer for $\Es$. 
Suppose now that $K$ is a rectifiable set in $\Om$ such that
$S(u) \subseteq K$ and 
\begin{equation}\label{gmzintro}
g^-(x,\nu_K(x))=0 \text{ for $\hn$-a.e. $x \in K$.}
\end{equation} 
Then we have immediately that the pair $(u,K)$ is a unilateral minimizer for $f$ and $g$ because for all pairs $(v,H)$ with $S(v) \subseteq H$ and $K \subseteq H$ we have
\begin{multline*}
\int_\Om f(x,\nabla u(x))\,dx=\Es(u) \le \Es(v)=
\int_\Om f(x,\nabla v(x))\,dx+\int_{S(v)}g^-(x,\nu)\,d\hn \\
=\int_\Om f(x,\nabla v(x))\,dx+\int_{S(v) \setminus K}g^-(x,\nu) 
\le \int_\Om f(x,\nabla v(x))\,dx+\int_{H \setminus K}g(x,\nu).
\end{multline*}
The rectifiable set $K$ satisfying \eqref{gmzintro} is provided in Section \ref{sigmaconv},
where we define a new variational notion of convergence for rectifiable sets which we call {\it $\sigma$-convergence}, and which departs from the
notion of $\sigma^p$-convergence given in \cite{DMFT}.
The $\sigma$-limit $K$ of a sequence of rectifiable sets $(K_n)_{n \in \N}$ 
is constructed looking for the $\Gamma$-limit $\hs^-$ in the strong topology of $L^1(\Om)$ of the functionals
\begin{equation*}
\label{hsnfuncintr}
\hs^-_n(u):=
\begin{cases}
\hn\left( S(u) \setminus K_n \right)  
& u \in P(\Om), \\
+\infty & \text{otherwise,}
\end{cases}
\end{equation*}
where $P(\Om)$ is the space of piecewise constant function in $\Om$ (see \eqref{pc}).
Roughly, the $\sigma$-limit $K$ is the maximal rectifiable set on which the density $h^-$ representing $\hs^-$ vanishes. By the growth estimate on $g$ it turns out that $K$ is also the maximal rectifiable set on which the density
$g^-$ vanishes, so that $K$ is the natural limit candidate for $K_n$ in order to preserve
the unilateral minimality property. The definition of $\sigma$-convergence involves only the surface energy densities $\hs^-_n$, and as
a consequence  it does not depend on the exponent $p$ and it is stable with respect to infinitesimal perturbations in length (see Remark \ref{kstrict}). Moreover 
it turns out that the $\sigma$-limit $K$ contains the $\sigma^p$-limit 
points of $(K_n)_{n \in \N}$, so that our $\Gamma$-convergence  approach improves also the minimality property given by the previous approaches. 
\par
Our method naturally extends to the case of varying bulk and surface energy densities $f_n$ and $g_n$, and this is indeed the main motivation for which we developed our $\Gamma$-convergence approach. The key point to recover effective energy densities $f$ and $g$ for the minimality property in the limit is a $\Gamma$-convergence result for functionals of the form
\begin{equation}
\label{freedisc}
\int_\Om f_n(x,\nabla u_n(x))\,dx+\int_{S(u_n)}g_n(x,\nu)\,d\hn(x).
\end{equation}
In Section \ref{gammasec}, we prove that the $\Gamma$-limit has the form
\begin{equation*}
\int_\Om f(x,\nabla u(x))\,dx+\int_{S(u)}g(x,\nu)\,d\hn(x),
\end{equation*}
where $f$ is determined only by $(f_n)_{n \in \N}$, and $g$ is determined only by $(g_n)_{n \in \N}$, that is no interaction occurs between the bulk and the surface part of the functionals in the $\Gamma$-convergence process. A result of this type has been proved in the case of periodic homogenization (in the vectorial case, and with dependence on the trace of $u$ in the surface part of the energy) by Braides, Defranceschi and Vitali \cite{BDV}. 
\par
We notice that an approach to stability in the line of Dal Maso, Francfort and Toader in the case of varying energies would have required a Transfer of Jump Sets for $f_n,g_n$ and $f,g$, which seems difficult to be derived directly. Our $\Gamma$-convergence approach also provides this result (Proposition \ref{transferofjump}).
\par
In section \ref{Composites} we deal with the study of quasistatic crack evolution in composite materials. More precisely we study the asymptotic behavior of a quasistatic evolution $t \to (u_n(t),K_n(t))$ relative to the bulk energy $f_n$ and the surface energy $g_n$. Using our stability result we prove (Theorem \ref{compevol}) that $t \to (u_n(t),K_n(t))$ converges to a quasistatic evolution $t \to (u(t),K(t))$ relative to 
the effective bulk and surface energy densities $f$ and $g$. Moreover convergence for bulk and surface energies for all times holds. This analysis applies to the case of composite materials, i.e. materials obtained through a fine mixture of different phases. The model case is that of periodic homogenization, i.e. materials with total energy given by
$$
\Es_\eps(u,K):=\int_\Om f \left( \frac{x}{\eps},\nabla u(x) \right)\,dx+
\int_K g \left( \frac{x}{\eps},\nu \right)\,d\hn(x),
$$
where $\eps$ is a small parameter giving the size of the mixture, and $f$, $g$ are periodic in $x$. Our result implies that a quasisistatic crack evolution $t \to (u_\eps(t),K_\eps(t))$ for $\eps$ small is very near to a quasistatic evolution for the homogeneous material having bulk and surface energy densities $f_{\rm hom}$ and $g_{\rm hom}$, which are obtained from $f$ and $g$ through  periodic homogenization formulas 
available in the literature (see for example \cite{BDV}).
\par
The paper is organized as follows. In Section \ref{prel} we make precise the functional setting of the problem. In Section \ref{blowupsec} we prove a blow up result for $\Gamma$-limits which will be employed in the proof of the main results. In Section \ref{reprsec} we prove some representation results which we use in Section \ref{gammasec} where we deal with the $\Gamma$-convergence of free discontinuity problems like \eqref{freedisc}. The notion of $\sigma$-convergence for rectifiable sets is contained in Section \ref{sigmaconv}, while the main result on stability for unilateral minimizers is contained in Section \ref{stabilitysec}. In Section \ref{stabilitysecbdry} we prove a stability result for unilateral minimality properties with boundary conditions which will be employed in Section \ref{Composites} for the study of quasistatic crack evolution in composite materials.

\section{The functional setting of the problem}
\label{prel}
We introduce now the precise functional setting for the study of the unilateral minimality property \eqref{minuni}. Throughout the paper we suppose that $\Om$ is a bounded open subset of $\R^N$ with Lipschitz boundary, and we denote by $\as(\Om)$ the family of its open subsets.
\par
In the unilateral minimality property \eqref{minuni}, we consider $(N-1)$-dimensional sets which are rectifiable, i.e. contained up to a set of $\hn$-measure zero in the union
of a sequence of $C^1$-hypersurfaces of $\R^N$. We will use the following notation: given $K_1,K_2$ rectifiable sets in $\R^N$, we say that 
$K_1 \tsub K_2$ if $K_1 \subseteq K_2$ up to a set of $\hs^{N-1}${-}measure zero; similarly we say that $K_1 \tilde{=} K_2$ if $K_1=K_2$
up to a set of $\hs^{N-1}${-}measure zero.
\par
Given $1<p<+\infty$, the functions 
in \eqref{minuni} belong to the space $SBV^p(\Om)$ defined as
$$
SBV^p(\Om):=\{u \in SBV(\Om)\,:\, \nabla u \in L^p(A,\R^N), \hn(S(u))<+\infty\}.
$$
For the notations and the general theory concerning the function space $SBV(\Om)$ ({\it special functions of bounded variation}), we refer the reader to \cite{AFP}. We will consider {\it weak convergence} in $SBV^p(\Om)$ defined in the following way: $u_n \weak u$ weakly in $SBV^p(\Om)$ if 
\begin{align*}
&u_n \to u \quad {strongly \; in}\; L^1(\Om), \\
&\nabla u_n \weak \nabla u \quad {weakly \; in}\;
L^p(\Om;\R^N), \\
&\hn(S(u_n)) \le C.
\end{align*}
\par
We indicate by $P(\Om)$ the family of sets with finite perimeter in $\Om$, that is 
the class of sets $E \subseteq \Om$ such that $1_E \in SBV(\Om)$. In view of the applications of Sections \ref{reprsec}, \ref{gammasec} and \ref{sigmaconv}, it will be useful to look at $P(\Om)$ in term of functions, that is to use the following equivalent description:
\begin{equation}
\label{pc}
P(\Om)=\{u \in SBV(\Om)\,:\, u(x) \in \{0,1\} \text{ for a.e. }x \in \Om\}.
\end{equation}

\section{Blow-up for $\Gamma$-limits}
\label{blowupsec}
In this section we state some blow-up results for $\Gamma$-convergent sequences
of integral functionals $\fs_n(u)$ which will be used in Section \ref{gammasec}. Moreover under additional hypothesis on $\fs_n$, we obtain a regularity result for the density of the $\Gamma$-limit $\fs$ which will be employed in Section \ref{Composites}. 
For the definition and the basic properties of $\Gamma$-convergence, we refer the reader to \cite{dm}.
\par
Let $1<p<+\infty$ and let $f: \Om \times \R^N \to [0,+\infty[$ be a Carath\'eodory function such that
\begin{equation}
\label{bulkenergyblow}
a_1(x)+\alpha |\xi|^p \le f(x,\xi) \le  a_2(x)+\beta |\xi|^p,
\end{equation}
where $a_1,a_2 \in L^1(\Om)$ and $\alpha,\beta>0$.
Let us assume that  
$$
\xi \to f(x,\xi) 
\quad
\text{is convex for a.e. }x \in \Om.
$$
Let $B_1$ be the unit ball in $\R^N$ with center $0$ and radius $1$. The following blow up result in the sense of $\Gamma$-convergence is a direct consequence of the Scorza-Dragoni theorem for Carath\'eodory functions and of \cite[Theorem 5.14]{dm}.

\begin{lemma}
\label{blowupthm}
Let $(\rho_k)_{k \in \N}$ be a sequence converging to zero. Then for a.e. $x \in \Om$ the functionals
\begin{equation*}
\label{frhon}
F_k(u):=
\begin{cases}
\int_{B_1} f(x+\rho_k y, \nabla u(y))\,dy      & u \in W^{1,p}(B_1), \\
+\infty      & \text{otherwise in }L^1(B_1)
\end{cases}
\end{equation*}
$\Gamma${-}converge in the strong topology of $L^1(B_1)$ to the functional
\begin{equation*}
\label{fx}
F(u):=
\begin{cases}
\int_{B_1} f(x, \nabla u(y))\,dy   & u \in W^{1,p}(B_1), \\
+\infty  & \text{otherwise in }L^1(B_1).
\end{cases}
\end{equation*}
\end{lemma}

Let us consider now $f_n: \Om \times \R^N \to [0,+\infty[$ Carath\'eodory function
satisfying the growth estimate \eqref{bulkenergyblow} uniformly in $n$, and
let $\fs_n:L^1(\Om) \times \as(\Om) \to [0,+\infty]$ be defined as
\begin{equation*}
\label{bulkfuncnblow}
\fs_n(u,A):=
\begin{cases}
\int_A f_n(x, \nabla u(x))\,dx      & u \in W^{1,p}(A), \\
+\infty      & \text{otherwise}.
\end{cases}
\end{equation*}
Let us assume (and this is always true up to a subsequence, see Theorem \ref{formulaf}) that for all $A \in \as(\Om)$ $\fs_n(\cdot,A)$
$\Gamma${-}converge with respect to the strong topology of $L^1(\Om)$
to a functional $\fs(\cdot,A)$ such that for all $u \in W^{1,p}(\Om)$
\begin{equation}
\label{bulkfuncblow}
\fs(u,A):=\int_A f(x, \nabla u(x))\,dx   
\end{equation}
for some Carath\'eodory function $f$ (independent of $u$ and $A$) which satisfies estimate \eqref{bulkenergyblow}.
Using Lemma \ref{blowupthm} and a diagonal argument we conclude that the following theorem holds.

\begin{theorem}
\label{blow2}
Let $(\rho_k)_{k \in \N}$ be a sequence converging to zero. Then for a.e. $x \in \Om$ there exists $(n_k)_{k \in \N}$ such that the functionals
\begin{equation*}
\label{frhon2}
F_k(u):=
\begin{cases}
\int_{B_1} f_{n_k}(x+\rho_k y, \nabla u(y))\,dy      & u \in W^{1,p}(B_1), \\
+\infty      & \text{otherwise in }L^1(B_1)
\end{cases}
\end{equation*}
$\Gamma${-}converge in the strong topology of $L^1(B_1)$ to the functional
\begin{equation*}
\label{fx2}
F(u):=
\begin{cases}
\int_{B_1} f(x, \nabla u(y))\,dy   & u \in W^{1,p}(B_1), \\
+\infty  & \text{otherwise in }L^1(B_1).
\end{cases}
\end{equation*}
\end{theorem}

\begin{remark}
\label{homrem}
{\rm
In the case of periodic homogenization, i.e. in the case in which
$f_n(x,\xi):=f(nx,\xi)$ with $f$ periodic in $x$, it is sufficient to choose $n_k$ in such a way that 
$n_k\rho_k \to +\infty$. In fact for $x=0$ we have
$$
F_k(u):=
\begin{cases}
\int_{B_1} f((n_k\rho_k)y, \nabla u(y))\,dy      & u \in W^{1,p}(B_1), \\
+\infty      & \text{otherwise in }L^1(B_1)
\end{cases}
$$
which still $\Gamma$-converges to (see for instance \cite{dm})
$$
F(u):=
\begin{cases}
\int_{B_1} f_{\rm hom}(\nabla u(y))\,dy      & u \in W^{1,p}(B_1), \\
+\infty      & \text{otherwise in }L^1(B_1).
\end{cases}
$$
}
\end{remark}

\vskip10pt
In the rest of the section we prove a regularity result for the density $f$ defined in \eqref{bulkfuncblow} under additional hypothesis on $f_n$ which will be employed in Section \ref{Composites}. 
Let us assume that for a.e. $x \in \Om$
\begin{itemize}
\item[(1)] $f_n(x,\cdot)$ is convex;
%\item[]
\vskip4pt
\item[(2)] $f_n(x,\cdot)$ is of class $C^1$;
\vskip4pt
\item[(3)] for all $M \ge 0$ and for all $\xi^1_n, \xi^2_n$ such that $|\xi^1_n| \le M$, $|\xi^1_n| \le M$, $|\xi^1_n-\xi^2_n| \to 0$ we have
\begin{equation}
\label{unifnablafxi}
|\nabla_\xi f_n(x,\xi^1_n)-\nabla_\xi f_n(x,\xi^2_n)| \to 0.
\end{equation}
\end{itemize}
Notice that for instance $f_n(x,\xi):=a_n(x)|\xi|^p$ with $\alpha \le a_n(x) \le \beta$ satisfies the assumptions above. Notice moreover that 
by lower semicontinuity of $\Gamma$-limits $\xi \to f(x,\xi)$ is convex for a.e. $x \in \Om$.
\par
We need the following lemma which is a straightforward variant of \cite[Lemma 4.9]{DMFT}.

\begin{lemma}
\label{convmom}
Let $(X, A, \mu)$ be a finite measure space, $p>1$,  $N \ge1$,
and let $H_n: X\times \R^N \to \R$ be a sequence of Carath\'eodory functions which satysfies the following properties: there exist 
a positive constant $a \ge 0$ and
a nonnegative function $b\in L^{p'}(X)$, with $p'=p/(p-1)$ such that
\begin{itemize}
\item[(1)] $|H_n(x,\xi)| \le a |\xi|^{p-1} + b(x)$ for every $x\in X, \, \xi\in\R^N$;
\vskip4pt
\item[(2)] for all $M \ge 0$ and for a.e. $x \in \Om$, for all $\xi^1_n, \xi^2_n$ such that $|\xi^1_n| \le M$, $|\xi^1_n| \le M$, $|\xi^1_n-\xi^2_n| \to 0$ we have
$$
|H_n(x,\xi^1_n)-H_n(x,\xi^2_n)| \to 0.
$$ 
\end{itemize}
Assume that $(\Phi_n)_{n \in \N}$ is bounded in $L^p(X,\R^N)$ and that $(\Psi_n)_{n \in \N}$ converges to $0$ strongly in $L^p(X,\R^N)$. Then
\begin{equation*}
\int_X [H_n(x,\Phi_n(x)+\Psi_n(x)) - H_n(x,\Phi_n(x))] \Phi(x) \, d\mu(x) \to 0,
\end{equation*}
for every $\Phi \in L^p(X,\R^N)$. 
\end{lemma}

The following regularity result on $f$ holds.

\begin{proposition}
\label{regolare} 
For a.e. $x\in \Om$ the function $\xi \to f(x,\xi)$ is of class $C^1$. 
 \end{proposition}

\begin{proof}
According to Theorem \ref{blow2}, let $x \in \Om$, $\rho_k \to 0$ and $(n_k)_{k \in \N}$ be such that $(F_k)_{k \in \N}$ $\Gamma${-}converges with respect to the strong topology of $L^1(B_1)$ to $F$.
\par
Let $(\phi_k)_{k \in \N}$ be a recovering sequence for the affine function $y \to \xi \cdot y$ with $\xi \in \R^N$. Up to a further subsequence, we can always assume that
there exists $\psi \in \R^N$ such that
\begin{equation}
\label{convpsi}
\frac{1}{|B_1|}\int_{B_1} \nabla_\xi f_{n_k}(x+\rho_k y,\nabla \phi_k(y))  \, dy \to \psi.
\end{equation}
Let $t_j \searrow 0$ and let $\eta \in \R^N$. By the convexity of $f_{n_k}$ in the second variable, we have
\begin{multline*}
\label{pp} 
\int_{B_1}f_{n_k} (x+\rho_k y, \nabla \phi_k(y) + t_j \eta)- f_{n_k} (x+\rho_k y, \nabla \phi_k(y)) \, dy \\ \le t_j\int_{B_1} \nabla_\xi f_{n_k}(x+\rho_k y, \nabla \phi_k(y)+t_j\eta)\eta \, dy. 
\end{multline*}
By $\Gamma$-convergence we can find $k_j$ such that
$$
\frac{f(x, \xi+ t_j \eta)- f(x,\xi)}{t_j} -\frac{1}{j} \le
\frac{1}{|B_1|}
\int_{B_1} \nabla_\xi f_{n_{k_j}}(x+\rho_{k_j}y, \nabla \phi_{k_j}(y)+t_j \eta)\eta\,dy,
$$
so that we have
\begin{equation}
\label{estfconvex}
\limsup_{j \to +\infty} \frac{f(x, \xi+ t_j \eta)- f(x,\xi)}{t_j}\le
\frac{1}{|B_1|}\limsup_{j \to +\infty}
\int_{B_1} \nabla_\xi f_{n_{k_j}}(x+\rho_{k_j}y, \nabla \phi_{k_j}(y)+t_j \eta)\eta\,dy.
\end{equation}
Notice that by Lemma \ref{convmom} and by \eqref{convpsi} we have that
\begin{multline*}
\lim_{j \to +\infty}
\int_{B_1} \nabla_\xi f_{n_{k_j}}(x+\rho_{k_j}y, \nabla \phi_{k_j}(y)+t_j \eta)\eta\,dy\\
=\lim_{j \to +\infty}
\int_{B_1} \nabla_\xi f_{n_{k_j}}(x+\rho_{k_j}y, \nabla \phi_{k_j}(y))\eta\,dy=
|B_1|\psi \eta,
\end{multline*}
and so for every subgradient $\zeta$ of $f(x,\cdot)$ at $\xi$ by \eqref{estfconvex} we have
$$
\zeta \eta  \le \limsup_{j \to +\infty} \frac{f(x, \xi+ t_j \eta)- f(x,\xi)}{t_j} \le
\psi\eta.
$$
We deduce that $\zeta=\psi$, so that $f(x,\cdot)$ is Gateaux differentiable at $\xi$ with
$\nabla_\xi f(x,\xi)=\psi$: since $f(x,\cdot)$ is convex, we get that $f(x,\cdot)$ is of class $C^1$.
\end{proof}

\begin{remark}
\label{noc1}
{\rm
An hypothesis of {\it equiuniform continuity} for $(\nabla_\xi f_n(x,\xi))_{n \in \N}$ like \eqref{unifnablafxi} is needed in order to preserve $C^1$-regularity in the passage from $f_n$ to $f$. Otherwise it is easy to provide a counterexample considering $\xi \to f_n(\xi)$ smooth convex functions uniformly converging to a non differentiable convex function $\xi \to f(\xi)$, and noting that the associated functionals $\Gamma$-converge.
}
\end{remark}

\section{Some integral representation lemmas}
\label{reprsec}
Let $a_1,a_2 \in L^1(\Om)$, $1<p<+\infty$, and let $\alpha,\beta>0$. For all $n \in \N$ let
$f_n: \Om \times \R^N \to [0,+\infty[$ be a Carath\'eodory function
such that for a.e. $x \in \Om$ and for all $\xi \in \R^N$
\begin{equation}
\label{bulkenergy}
a_1(x)+\alpha |\xi|^p \le f_n(x,\xi) \le  a_2(x)+\beta |\xi|^p,
\end{equation}
and let $g_n: \Om \times S^{N-1} \to [0,+\infty[$ be a Borel function
such that for $\hn$-a.e. $x \in \Om$ and for all $\nu \in S^{N-1}:=\{\eta \in \R^N:|\eta|=1\}$
\begin{equation}
\label{surfenergy}
\alpha \le g_n(x,\nu) \le \beta.
\end{equation}
\par
In Section \ref{gammasec} we will be interested in the functionals on $L^1(\Om) \times \as(\Om)$
\begin{equation}
\label{funcn2}
\Es_n(u,A):=
\begin{cases}
\int_A f_n(x,\nabla u(x))\,dx+ 
\int_{A \cap (S(u) \setminus K_n)} g_n(x, \nu)\,d\hn(x)  
& u \in SBV^p(A), \\
+\infty & \text{otherwise},
\end{cases}
\end{equation}
where $\as(\Om)$ denotes the family of open subsets of $\Om$, and $(K_n)_{n \in \N}$ is a sequence of rectifiable sets in $\Om$ such that
\begin{equation*}
\label{boundkn}
\hn(K_n) \le C.
\end{equation*}
In particular we will be interested in the $\Gamma${-}limit in the strong topology of $L^1(\Om)$ of $(\Es_n(\cdot,A))_{n \in \N}$ for every $A \in \as(\Om)$.
To this extend we consider the functionals $\fs_n:L^1(\Om) \times \as(\Om) \to [0,+\infty]$
\begin{equation}
\label{bulkfuncn}
\fs_n(u,A):=
\begin{cases}
\int_A f_n(x, \nabla u(x))\,dx      & u \in W^{1,p}(A), \\
+\infty      & \text{otherwise},
\end{cases}
\end{equation}
and the functionals $\gs^-_n:P(\Om) \times \as(\Om) \to [0,+\infty[$
\begin{equation}
\label{surffuncnm}
\gs^-_n(u,A):=
\int_{A \cap (S(u) \setminus K_n)} g_n(x, \nu)\,d\hn(x)  
\end{equation}
defined on Sobolev and piecewise constant functions with values in $\{0,1\}$ (see \eqref{pc}) respectively, and we will reconstruct the $\Gamma${-}limit of $(\Es_n(\cdot,A))_{n \in \N}$ through the $\Gamma${-}limits of $(\fs_n(\cdot,A))_{n \in \N}$ and $(\gs^-_n(\cdot,A))_{n \in \N}$. 
\par
For the results of Section \ref{stabilitysec}, we will need also the functionals
$\gs_n:P(\Om) \times \as(\Om) \to [0,+\infty[$
\begin{equation}
\label{surffuncn}
\gs_n(u,A):=
\int_{A \cap S(u)} g_n(x, \nu)\,d\hn(x).  
\end{equation}
\par
In the following, for every functional $\hs$ defined on $X \times \as(\Om)$ with $X=L^1(\Om)$ or $X=P(\Om)$ with values in $[0,+\infty]$, for every $\psi \in L^1(A)$ and $A \in \as(\Om)$ we will use the notation
\begin{equation}
\label{mbf}
{\bf m}_{\hs}(\psi,A)=\inf_{u \in X} \{\hs(u,A)\,:\,u=\psi \text{ in a neighborhood of }\partial A\}.
\end{equation}
Moreover for all $x \in \R^N$, $a,b \in \R$ and $\nu \in S^{N-1}$  let $u_{x,a,b,\nu}:B_1(x) \to \R$ be defined by
\begin{equation}
\label{uxabnu}
u_{x,a,b,\nu}(y):=
\begin{cases}
b & \text{if }(y-x)\nu \ge 0, \\
a & \text{if }(y-x)\nu < 0,
\end{cases}
\end{equation}
where $B_1(x)$ is the ball of center $x$ and radius $1$.
\par
The following $\Gamma$-convergence and representation result for the functionals $\fs_n$ holds (see Buttazzo and Dal Maso \cite{BDM}, Bouchitt\'e, Fonseca, Leoni and Mascarenhas \cite[Theorem 2]{BFLM}).

\begin{proposition}
\label{formulaf}
There exists $\fs:L^1(\Om) \times \as(\Om) \to [0,+\infty]$ such that up to a subsequence the functionals $\fs_n(\cdot,A)$ $\Gamma$-converge in the strong topology of $L^1(\Om)$ to $\fs(\cdot,A)$ for every $A \in \as(\Om)$. Moreover for all $u \in W^{1,p}(\Om)$ we have that
\begin{equation}
\label{bulkfunc}
\fs(u,A)= \int_A f(x, \nabla u(x))\,dx,
\end{equation}
where
\begin{equation}
\label{reprf}
f(x,\xi):= \limsup_{\rho \to 0^+}
\frac{{\bf m}_{\fs}(\xi(z-x),B_\rho(x))}{\omega_N \rho^N},
\end{equation}
$\bf{m}_\fs$ is defined in \eqref{mbf}, and $\omega_N$ is the volume of the unit ball in $\R^N$. Finally $f$ is a Carath\'eodory function satisfying the growth conditions \eqref{bulkenergy}.
\end{proposition}

Let us come to the functionals $\gs_n$ defined in \eqref{surffuncn}. The following proposition holds (see Ambrosio and Braides \cite[Theorem 3.2]{AB1}, Bouchitt\'e, Fonseca, Leoni and Mascarenhas \cite[Theorem 3]{BFLM}).

\begin{proposition}
\label{formulag2}
There exists $\gs: P(\Om) \times \as(\Om) \to [0,+\infty[$ such that up to a subsequence $\gs_n(\cdot,A)$ $\Gamma$-converge with respect to the strong topology of $L^1(\Om)$ to $\gs(\cdot,A)$ for all $A \in \as(\Om)$. Moreover for all $u \in P(\Om)$ and $A \in \as(\Om)$ we have that
\begin{equation}
\label{surffunc2}
\gs(u,A)=\int_{A \cap S(u)} g(x, \nu)\,dx,      
\end{equation}
with
\begin{equation}
\label{reprg2}
g(x,\nu):= \limsup_{\rho \to 0^+}
\frac{{\bf m}_{\gs}(u_{x,0,1,\nu},B_\rho(x))}{\omega_{N-1} \rho^{N-1}},
\end{equation}
where ${\bf m}_{\gs}$ is defined in \eqref{mbf} and $u_{x,0,1,\nu}$ is as in \eqref{uxabnu}.
\end{proposition}

Let us come to the functionals $\gs^-_n$ defined in \eqref{surffuncnm}. The following proposition holds.

\begin{proposition}
\label{formulag}
There exists $\gs^-: P(\Om) \times \as(\Om) \to [0,+\infty[$ such that up to a subsequence $\gs^-_n(\cdot,A)$ $\Gamma$-converge with respect to the strong topology of $L^1(\Om)$ to $\gs^-(\cdot,A)$ for all $A \in \as(\Om)$. Moreover for all $u \in P(\Om)$ and $A \in \as(\Om)$ we have that
\begin{equation}
\label{surffunc}
\gs^-(u,A)=
\int_{A \cap S(u)} g^-(x, \nu)\,d\hn(x),  
\end{equation}
with
\begin{equation}
\label{reprgm}
g^-(x,\nu):= \limsup_{\rho \to 0^+}
\frac{{\bf m}_{\gs^-}(u_{x,0,1,\nu},B_\rho(x))}{\omega_{N-1} \rho^{N-1}},
\end{equation}
where ${\bf m}_{\gs^-}$ is defined in \eqref{mbf} and $u_{x,0,1,\nu}$ is as in \eqref{uxabnu}.
\end{proposition}

\begin{proof}
The $\Gamma$-convergence result for $\gs_n^-(\cdot,A)$ is given by the result of Ambrosio and Braides \cite{AB1}. For the sequel we need also the explicit formula \eqref{reprgm} for the density $g^-$ which is not given directly by the results of \cite{AB1} and \cite{BFLM} because of a lack of coercivity from below. Let us briefly sketch how to prove that $g^-$ defined in \eqref{reprgm} represents $\gs^-$.
According to Proposition \ref{formulag2}, let us consider the densities $g^\eps(x,\nu)$ representing the $\Gamma$-limit 
$\gs^\eps(\cdot,A)$ of the (uniformly coercive) functionals
\begin{equation*}
\label{surffuncneps}
\gs^\eps_n(u,A):=
\int_{A \cap S(u)} g^\eps_n(x, \nu)\,d\hn(x),  
\end{equation*}
where
\begin{equation*}
\label{defgneps}
g^\eps_n(x,\nu):=
\begin{cases}
\eps  & \text{if } x \in K_n, \nu=\nu_{K_n}(x), \\
g_n(x,\nu)  & \text{otherwise}.
\end{cases}
\end{equation*}
We have immediately that $\gs^\eps(u,A) \to \gs^-(u,A)$ as $\eps \to 0$ for all $u \in P(\Om)$ and $A \in \as(\Om)$. Let $\mu$ be the \wstar limit of $\hn \res K_n$ (up to a subsequence) in the sense of measures. Notice that (see for instance \cite[Theorem 2.56]{AFP})
up to a set of $\hn${-}measure zero we have 
\begin{equation*}
\label{hx}
H(x):=\limsup_{\rho \to 0^+} \frac{\mu(\bar B_\rho(x))}{\omega_{N-1} \rho^{N-1}}
<+\infty.
\end{equation*}
Then the result follows noting that for all $x \in \Om$ with $H(x)<+\infty$ we have 
$$
g^-(x,\nu)=\lim_{\eps \to 0}g^\eps(x,\nu).
$$
\end{proof}

\begin{remark}
\label{gnmpc}
{\rm
It is immediate to check that if we replace $P(\Om)$ in Proposition \ref{formulag} by the space $P_{a,b}(\Om):=\{u \in SBV(\Om)\,:\,u(x) \in \{a,b\} \text{ for a.e. }x \in \Om\}$, with
$a,b \in \R$, then the $\Gamma$-limit in the strong topology of $L^1(\Om)$
of $\gs^-_n(\cdot,A)$ can still be represented by the density $g^-$ defined in \eqref{reprgm}.
}
\end{remark}

Let us finally come to the functionals $\Es_n$ defined in \eqref{funcn2}. Using the growth estimates \eqref{bulkenergy} and \eqref{surfenergy} on $f_n$ and $g_n$ (see \cite{BDV}), there exists $\Es: L^1(\Om) \times \as(\Om) \to [0,+\infty[$ such that up to a subsequence $\Es_n(\cdot,A)$ $\Gamma$-converge in the strong topology of $L^1(\Om)$ to $\Es(\cdot,A)$ for all $A \in \as(\Om)$. 
For every $\eps>0$ let us set
$$
\Es_\eps(u,A):=\Es(u,A)+\eps \int_{S(u) \cap A} 1+|[u]| \,d\hn,
$$
where $[u](x)$ denotes the jump of $u$ at $x$, i.e. $[u](x):=u^+(x)-u^-(x)$.
By the representation result of Bouchitt\'e, Fonseca, Leoni and Mascarenhas \cite[Theorem 1]{BFLM} we get that for all $u \in SBV^p(\Om)$ and $A \in \as(\Om)$
$$
\Es_\eps(u,A)=
\int_A f^\eps_\infty(x,\nabla u(x))\,dx+ \int_{A \cap S(u)} g^\eps_\infty(x,u^-(x),u^+(x),\nu)\,d\hn(x)
$$
with $f^\eps_\infty$ and $g^\eps_\infty$ satisfying the following formulas
\begin{equation}
\label{reprfinfty}
f^\eps_\infty(x,\xi):= \limsup_{\rho \to 0^+}
\frac{{\bf m}_{\Es_\eps}(\xi(z-x),B_\rho(x))}{\omega_N \rho^N},
\end{equation}
\begin{equation}
\label{reprginfty}
g^\eps_\infty(x,a,b,\nu):= \limsup_{\rho \to 0^+}
\frac{{\bf m}_{\Es_\eps}(u_{x,a,b,\nu},B_\rho(x))}{\omega_{N-1} \rho^{N-1}},
\end{equation}
where ${\bf m}_{\Es_\eps}$ is defined in \eqref{mbf} and $u_{x,a,b,\nu}$ is as in \eqref{uxabnu}.
\par
Notice that $f^\eps_\infty$ and $g^\eps_\infty$ are monotone decreasing in $\eps$, and that $\Es_\eps(\cdot,A)$ converges pointwise to $\Es(\cdot,A)$ as $\eps \to 0$ for all $A \in \as(\Om)$. We conclude that the representation result for $\Es_\eps$ implies a representation result for the functional $\Es$.
\par
Summarizing we have that the following proposition holds.

\begin{proposition}
\label{formulafginfty}
There exists $\Es:L^1(\Om) \times \as(\Om) \to [0,+\infty]$ such that up to a subsequence $\Es_n(\cdot,A)$ $\Gamma$-converges in the strong topology of $L^1(\Om)$ to $\Es(\cdot,A)$ for every $A \in \as(\Om)$. Moreover,
for every $u \in SBV^p(\Om)$ and $A \in \as(\Om)$ we have that
\begin{equation*}
\label{esrepr}
\Es(u,A)=
\int_A f_\infty(x,\nabla u(x))\,dx+ \int_{A \cap S(u)} g_\infty(x,u^-(x),u^+(x),\nu)\,d\hn(x)
\end{equation*}
with
\begin{equation}
\label{fepstof}
f_\infty(x,\xi):=\lim_{\eps \to 0}f^\eps_\infty(x,\xi)
\qquad\text{and}\qquad
g_\infty(x,a,b,\nu):=\lim_{\eps \to 0}g^\eps_\infty(x,a,b,\nu),
\end{equation}
where $f^\eps_\infty$ and $g^\eps_\infty$ are defined in \eqref{reprfinfty} and \eqref{reprginfty} respectively.
\end{proposition}

\begin{remark}
\label{remlocalization}
{\rm
In the rest of the paper we will often make use the following property which is implied by the fact that $\Es(u,\cdot)$ is a Radon measure for every $u \in SBV^p(\Om)$. If $(u_n)_{n \in \N}$ is a recovering sequence for $u$ with respect to $\Es_n(\cdot, \Om)$, then $(u_n)_{n \in \N}$ is optimal for $u$ with respect to $\Es_n(\cdot,A)$ for every $A \in \as(\Om)$ such that the measure $\Es(u,\cdot)$ vanishes on $\partial A$.
}
\end{remark}

\section{A $\Gamma${-}convergence result for free discontinuity problems}
\label{gammasec}
The main result of this section is the following $\Gamma$-convergence theorem concerning the functionals $\Es_n$ defined in \eqref{funcn2}.

\begin{theorem}
\label{gammaresult}
Let $(K_n)_{n \in \N}$ be a sequence of rectifiable sets in $\Om$ such that $\hn(K_n) \le C$ for all $n \in \N$.
Let us assume that for all $A \in \as(\Om)$ the functionals 
$\fs_n(\cdot,A)$ and $\gs^-_n(\cdot,A)$ defined in \eqref{bulkfuncn} and \eqref{surffuncnm} $\Gamma${-}converge in the strong topology of $L^1(\Om)$ to $\fs(\cdot,A)$ and $\gs^-(\cdot,A)$ respectively.
Then for all $A \in \as(\Om)$ the functionals $\Es_n(\cdot,A)$ defined in \eqref{funcn2} $\Gamma$-converge in the strong topology of $L^1(\Om)$ to $\Es(\cdot,A)$ 
such that for all $u \in SBV^p(\Om)$ and $A \in \as(\Om)$
\begin{equation*}
\label{func}
\Es(u,A)=
\int_A f(x,\nabla u(x))\,dx+ 
\int_{A \cap S(u)} g^-(x, \nu)\,d\hn(x),  
\end{equation*}
where $f$ and $g^-$ are the densities of $\fs$ and $\gs^-$ according to Propositions \ref{formulaf} and \ref{formulag}.
\end{theorem}

\begin{proof}
We know that up to a subsequence the functionals $\Es_n(\cdot,A)$ $\Gamma$-converge in the strong topology of $L^1(\Om)$ to a functional $\Es(\cdot,A)$ for every $A \in \as(\Om)$, and that by Proposition \ref{formulafginfty} for all $u \in SBV^p(\Om)$ and for all $A \in \as(\Om)$ we have
$$
\Es(u,A)=\int_A f_\infty(x,\nabla u)\,dx+
\int_{S(u) \cap A}g_\infty(x,u^-(x),u^+(x),\nu)\,d\hn(x),
$$
where $f_\infty$ and $g_\infty$ satisfy formula \eqref{fepstof}.
The theorem will be proved if we show that for all $u \in SBV^p(\Om)$
we have $f_\infty(x, \nabla u(x))=f(x, \nabla u(x))$ for a.e. $x \in \Om$, and
$g_\infty(x,u^-(x),u^+(x),\nu_{S(u)}(x))=g^-(x,\nu_{S(u)}(x))$
for $\hn${-}a.e. $x \in S(u)$, where $\nu_{S(u)}(x)$ is the normal to $S(u)$ at $x$.
\par
The proof will be divided into four steps.

\vskip10pt\par\noindent
{\bf Step 1: 
$\mathbf{f_\infty(x, \nabla u(x)) \le f(x, \nabla u(x))}$ for a.e. $\mathbf{x \in \Om}$.}
\par
This inequalty can be derived using the explicite formulas for $f_\infty$ and
$f$. Let $x \in \Om$, $\xi \in \R^N$, and let us fix $\eps>0$. For every $\rho>0$ let $u_{\eps,\rho} \in W^{1,p}(B_\rho(x))$ be such that $u_{\eps,\rho}(z)=\xi(z-x)$ in a neighborhood of $\partial B_\rho(x)$ and
$$
\fs(u_{\eps,\rho},B_\rho(x)) \le {\bf m}_\fs(\xi(z-x),B_\rho(x))+\eps \omega_N \rho^N.
$$
Then we get
\begin{multline*}
f^\eps_\infty(x,\xi)= \limsup_{\rho \to 0^+} 
\frac{{\bf m}_{\Es_\eps}(\xi(z-x),B_\rho(x))}{\omega_N \rho^N} 
\le \limsup_{\rho \to 0^+}
\frac{\Es(u_{\eps,\rho},B_\rho(x))}{\omega_N \rho^N} \\
\le \limsup_{\rho \to 0^+}
\frac{\fs(u_{\eps,\rho},B_\rho(x))}{\omega_N \rho^N} 
\le \limsup_{\rho \to 0^+}
\frac{{\bf m}_\fs(\xi(z-x),B_\rho(x))}{\omega_N \rho^N} +\eps=f(x,\xi) +\eps.
\end{multline*}
Letting $\eps \to 0$, we obtain that $f_\infty(x,\xi) \le f(x,\xi)$, so that the step is concluded.

\vskip10pt\par\noindent
{\bf Step 2: 
$\mathbf{f_\infty(x, \nabla u(x)) \ge f(x, \nabla u(x))}$ for a.e. $\mathbf{x \in \Om}$.}
\par
We can consider those $x \in \Om$ such that $u$ is approximatively differentiable at $x$,
$x$ is a Lebesgue point for $f(\cdot,\xi)$ for all $\xi \in \R^N$ and such that
\begin{equation}
\label{finfty}
f_\infty(x, \nabla u(x))= 
\lim_{\rho \to 0^+} \frac{\Es(u,B_\rho(x))}{\omega_N \rho^N}<+\infty.
\end{equation}
Let moreover $(u_n)_{n \in \N}$ be a recovering sequence for $\Es(u,\Om)$: by \eqref{surfenergy} and since $\hn(K_n) \le C$, we have that
$\hn(S(u_n))$ is bounded and so up to a subsequence
$$
\mu_n:=\hn \res S(u_n) \weakst \mu
\qquad
\text{\wlystar in the sense of measures}
$$
for some Borel measure $\mu$. We can assume that (see for instance \cite[Theorem 2.56]{AFP})
\begin{equation}
\label{densitymu}
\limsup_{\rho \to 0^+} \frac{\mu(\bar B_\rho(x))}{\rho^{N-1}}=0.
\end{equation}
Let $\rho_i \searrow 0$ be such that $\Es(u,\partial B_{\rho_i}(x))=0$. In view of Remark \ref{remlocalization}, for every $i$ there exists $n_i$ such that for 
$n \ge n_i$
\begin{multline*}
\label{ineq1}
\frac{\Es(u,B_{\rho_i}(x))}{\omega_N \rho_i^N}
\ge 
\frac{\Es_n(u_n,B_{\rho_i}(x))}{\omega_N \rho_i^N}
-\frac{1}{i} \\
\ge \frac{\int_{B_{\rho_i}(x)} f_n(x, \nabla u_n(x))\,dx}{\omega_N \rho_i^N}
-\frac{1}{i} 
=\frac{1}{\omega_N}\int_{B_1} f_n(x+\rho_i y, \nabla v^i_n(y))\,dy
-\frac{1}{i} 
\end{multline*}
where
$$
v^i_n(y):= \frac{u_n(x+\rho_i y)-u(x)}{\rho_i}.
$$
Taking into account the assumptions on $x$, \eqref{finfty} and \eqref{densitymu}, we can choose $(n_i)_{i \in \N}$ is such a way that
\begin{equation*}
\label{vinl1}
v^i_{n_i} \to \nabla u(x) \cdot y
\quad\quad
\text{strongly in }L^1(B_1) \text{ for }i \to +\infty,
\end{equation*}
\begin{equation*}
\label{nablavin}
(\nabla v^i_{n_i})_{i \in \N} \text{ is bounded in }L^p(B_1,\R^N),
\end{equation*}
\begin{equation*}
\label{hvinzero}
\lim_{i \to +\infty}\hn(S(v^i_{n_i}))=0,
\end{equation*}
and 
\begin{equation}
\label{energyvin}
f_{\infty}(x,\nabla u(x))=\lim_{i \to +\infty}\frac{\Es(u,B_{\rho_i}(x))}{\omega_N \rho_i^N} \ge
\liminf_{i \to +\infty}\frac{1}{\omega_N}\int_{B_1} f_{n_i}(x+\rho_i y, \nabla v^i_{n_i}(y))\,dy.
\end{equation}
Moreover by a truncation argument we can assume that $(v^i_{n_i})_{i \in \N}$ is uniformly bounded in $L^\infty(B_1)$, 
so that we get
$$
\|\nabla v^i_{n_i}\|^p_{L^p(B_1,\R^N)}+ \int_{S(v^i_{n_i})} |[v^i_{n_i}]| \,d\hn \le C
\quad\text{ and }\quad
\lim_{i \to +\infty}\hn(S(v^i_{n_i}))=0.
$$
Following Kristensen \cite{K} we get that there exists $w_i \in W^{1,\infty}(B_1)$ such that $w_i \to \nabla u(x)\cdot y$ strongly in $L^1(B_1)$ as $i \to +\infty$ and such that
\begin{equation*}
\label{convsob}
\liminf_{i \to +\infty} \int_{B_1} f_{n_i}(x+\rho_i y, \nabla v^i_{n_i}(y))\,dy
=
\liminf_{i \to +\infty} \int_{B_1} f_{n_i}(x+\rho_i y, \nabla w_i(y))\,dy.
\end{equation*}
If $n_i$ is choosen such that the blow-up for $\Gamma$-limits given by Theorem \ref{blow2} holds, we get that
$$
\liminf_{i \to +\infty} \int_{B_1} f_{n_i}(x+\rho_i y, \nabla w_i(y))\,dy
\ge \omega_N f(x,\nabla u(x)),
$$
so that in view of \eqref{energyvin} we obtain
$$
f_\infty(x,\nabla u(x)) \ge f(x, \nabla u(x)).
$$

\vskip10pt\par\noindent
{\bf Step 3: 
$\mathbf{g_\infty(x,u^-(x),u^+(x),\nu_{S(u)}(x)) \le g^-(x, \nu_{S(u)}(x))}$ for $\mathbf{\hn}${-}a.e. $\mathbf{x \in S(u)}$.}
\par
Up to a subsequence, we have that
$$
\mu_n:=\hn \res K_n \weakst \mu
$$
\wlystar in the sense of measures. Since $\hn(K_n) \le C$ we have that for
$\hn$-a.e. $x \in \Om$ (see for instance \cite[Theorem 2.56]{AFP})
\begin{equation}
\label{hfinite}
H(x):=\limsup_{\rho \to 0^+} \frac{\mu(\bar B_\rho(x))}{\omega_{N-1} \rho^{N-1}}<+\infty.
\end{equation}
We claim that for all $v \in P(\Om)$ and $A \in \as(\Om)$ such that $\bar A \subseteq \Om$ 
\begin{equation}
\label{comp1}
\alpha \hn(S(v) \cap A) \le \gs^-(v,A)+\alpha \mu(\bar A).
\end{equation}
In fact we have that for all $n \in \N$
$$
\alpha \hn\left( (S(v) \setminus K_n) \cap A \right) \le \gs^-_n(v,A)
$$
so that
$$
\alpha \hn\left( S(v) \cap A \right) \le \gs^-_n(v,A)+\alpha \mu_n(A)
$$
and so passing to the $\Gamma$-limit for $n \to +\infty$ we obtain that \eqref{comp1} holds.
\par
Let us choose $x \in S(u)$ in such a way that \eqref{hfinite} holds and such that
$$
\limsup_{\rho \to 0^+} \frac{\int_{B_\rho(x)}a_2\,dx}{\rho^{N-1}}=0,
$$
where $a_2$ is defined in \eqref{bulkenergy}.
Let us indicate $u^-(x), u^+(x)$ and $\nu_{S(u)}(x)$ simply by $u^-, u^+$ and $\nu$. Let us moreover set $[u]:=u^+-u^-$.
\par
Following Remark \ref{gnmpc}, let us consider the functionals $\gs^-_n$ defined in \eqref{surffuncnm} acting on the space $P_{u^-,u^+}(\Om):=\{u \in SBV(\Om)\,:\, u(y) \in \{u^-,u^+\} \text{ for a.e. }y \in \Om\}$.
\par
Let us fix $\eps>0$. For every $\rho>0$, let $u_{\eps,\rho} \in P_{u^-,u^+}(B_\rho(x))$ be such that $u_{\eps,\rho}=u_{x,u^-,u^+,\nu}$ in a neighborhood of $B_\rho(x)$ and 
$$
\gs^-(u_{\eps,\rho},B_\rho(x)) \le {\bf m}_{\gs^-}(u_{x,u^-,u^+,\nu},B_\rho(x))+\eps \omega_{N-1} \rho^{N-1}.
$$
Then we get in view of \eqref{reprginfty} and \eqref{comp1}
\begin{multline*}
g^\eps_\infty(x,u^-,u^+,\nu)=
\limsup_{\rho \to 0^+}
\frac{{\bf m}_{\Es_\eps}(u_{x,u^-,u^+,\nu},B_\rho(x))}{\omega_{N-1} \rho^{N-1}} \\
\le
\limsup_{\rho \to 0^+}
\frac{\Es(u_{\eps,\rho},B_\rho(x))+\eps(1+|[u]|)\hn(S(u_{\eps,\rho}) \cap B_\rho(x))}{\omega_{N-1} \rho^{N-1}} \\
\le \limsup_{\rho \to 0^+}
\frac{\int_{B_\rho(x)}a_2\,dx+\gs^-(u_{\eps,\rho},B_\rho(x))+\frac{\eps}{\alpha}(1+|[u]|)
(\gs^-(u_{\eps,\rho},B_\rho(x))+\alpha \mu(\bar B_\rho(x)))}{\omega_{N-1} \rho^{N-1}} \\
\le
\limsup_{\rho \to 0^+}
\frac{(1+\frac{\eps}{\alpha}+\frac{\eps}{\alpha}|[u]|) ({\bf m}_{\gs^-}
(u_{x,u^-,u^+,\nu},B_\rho(x))+\eps \omega_{N-1}\rho^{N-1})
+\eps(1+|[u]|)\mu(\bar B_\rho(x)))}{\omega_{N-1} \rho^{N-1}} \\
\le 
\left( 1+\frac{\eps}{\alpha}+\frac{\eps}{\alpha}|[u]| \right)
(g^-(x,\nu)+\eps)+\eps(1+|[u]|)H(x).
\end{multline*}
Letting $\eps \to 0$ we obtain $g_\infty(x,u^-,u^+,\nu) \le g^-(x,\nu)$, so that the step is concluded.

\vskip10pt\par\noindent
{\bf Step 4: 
$\mathbf{g_\infty(x,u^-(x),u^+(x),\nu_{S(u)}(x)) \ge g^-(x, \nu_{S(u)}(x))}$ for $\mathbf{\hn}${-}a.e. $\mathbf{x \in S(u)}$.}
\par
Let us choose $x \in S(u)$ which is an approximate jump point for $u$,
\begin{equation}
\label{ginftyux}
g_\infty(x,u^-(x),u^+(x),\nu_{S(u)}(x))=\lim_{\rho \to 0^+} 
\frac{\Es(u,B_\rho(x))}{\omega_{N-1}\rho^{N-1}}<+\infty,
\end{equation}
and such that
\begin{equation}
\label{a1rhoi}
\lim_{\rho \to 0^+} \frac{\int_{B_\rho(x)}|a_1(y)|\,dy}{\rho^{N-1}}=0,
\end{equation}
where $a_1$ is defined in \eqref{bulkenergy}.
\par
Since $\hn(K_n) \le C$, up to a subsequence we have
$$
\mu_n:= \hn \res K_n \weakst \mu
\qquad
\text{\wlystar in the sense of measures}
$$ 
for some Borel measure $\mu$.
We can assume that (see for instance \cite[Theorem 2.56]{AFP})
\begin{equation}
\label{muxfinite}
\limsup_{\rho \to 0^+} \frac{\mu(B_\rho(x))}{\rho^{N-1}}<+\infty.
\end{equation}
Let $(u_n)_{n \in \N}$ be a recovering sequence for $\Es(u,\Om)$, and let $\rho_i \searrow 0$ be such that $\Es(u,\partial B_{\rho_i}(x))=0$. For every $i \in \N$ there exists $n_i \in \N$ such that for $n \ge n_i$ we have
\begin{multline}
\label{ineqstep3}
\frac{\Es(u,B_{\rho_i}(x))}{\omega_{N-1}\rho_i^{N-1}} \ge 
\frac{\Es_n(u_n,B_{\rho_i}(x))}{\omega_{N-1}\rho_i^{N-1}} -\frac{1}{i} \\
\ge 
\frac{\int_{B_{\rho_i}(x) \cap [S(u_n) \setminus K_n]}g_n(x,\nu)\,d\hn(x)}
{\omega_{N-1}\rho_i^{N-1}} 
+\frac{\int_{B_{\rho_i}(x)}a_1(y)\,dy}{\omega_{N-1}\rho_i^{N-1}}
-\frac{1}{i} \\
=
\frac{1}{\omega_{N-1}} \int_{B_1 \cap [S(v^i_n) \setminus K^i_n]}g_n(x+\rho_iy,\nu)\,d\hn(y)
+\frac{\int_{B_{\rho_i}(x)}a_1(y)\,dy}{\omega_{N-1}\rho_i^{N-1}}
 -\frac{1}{i},
\end{multline}
where
\begin{equation*}
\label{vin}
v^i_n(y):=u_n(x+\rho_i y)
\quad
\text{ and }
\quad
K^i_n:= \frac{\{K_n \cap B_{\rho_i}(x)\}-x}{\rho_i}.
\end{equation*}
We claim that we can find $w^i_n$ piecewise constant in $B_1$ such that
for $n \to +\infty$
$$
w^i_n \to w^i
\quad\quad
\text{strongly in }L^1(B_1),
$$
where $w^i$ is piecewise constant and $w^i=u_{0,0,1,\nu_{S(u)}(x)}$ in a neighborhood of the boundary, and such that for $n$ large 
\begin{equation*}
\label{ineqstep3tris}
\int_{B_1 \cap [S(v^i_n) \setminus K^i_n]}g_n(x+\rho_iy,\nu)\,d\hn(y)
\ge
\int_{B_1 \cap [S(w^i_n) \setminus K^i_n]}g_n(x+\rho_iy,\nu)\,d\hn(y)
-e_i,
\end{equation*}
with $e_i \to 0$ for $i \to +\infty$.
\par
Using the claim, by \eqref{ginftyux}, \eqref{ineqstep3}, \eqref{ineqstep3tris} and \eqref{a1rhoi} we have that for $n$ large
$$
g_\infty(x,u^-(x),u^+(x),\nu_{S(u)}(x)) \ge 
\frac{\int_{B_{\rho_i} \cap [S(z^i_n) \setminus K_n]}g_n(\zeta,\nu)\,d\hn(\zeta)}{\omega_{N-1}\rho_i^{N-1}}-\hat e_i=
\frac{\gs^-_n(z^i_n,B_{\rho_i}(x))}{\omega_{N-1}\rho_i^{N-1}}-\hat e_i
$$
where $\hat e_i \to 0$ and
$$
z^i_n(\zeta):=w^i_n \left( \frac{\zeta-x}{\rho_i} \right)
\to z^i(\zeta):=w^i \left( \frac{\zeta-x}{\rho_i} \right)
\quad\quad
\text{strongly in }L^1(B_{\rho_i}(x)).
$$
By the $\Gamma${-}convergence assumption on $\gs^-_n$, using $\Gamma$-liminf inequality we have that
$$
g_\infty(x,u^-(x),u^+(x),\nu_{S(u)}(x)) \ge 
\frac{\gs^-(z^i,B_{\rho_i}(x))}
{\omega_{N-1}\rho_i^{N-1}}-\hat e_i
\ge
\frac{{\bf m}_{\gs^-}(u_{x,0,1,\nu_{S(u)}(x)},B_{\rho_i})}
{\omega_{N-1}\rho_i^{N-1}}-\hat e_i.
$$
Letting $i \to +\infty$, and recalling the representation formula \eqref{reprgm} for $g^-(x,\nu)$, we have that the result is proved.
\par
In order to complete the proof of the step, we have to prove the claim.
Since
$$
\nabla v^i_n(y)=\rho_i \nabla u_n(x+\rho_i y),
$$
we get by the coercivity assumption \eqref{bulkenergy}
\begin{multline*}
\int_{B_1}|\nabla v^i_n(y)|^p\,dy=
\rho_i^p \int_{B_1}|\nabla u_n(x+\rho_i y)|^p\,dy 
= \rho_i^{p} \frac{\int_{B_{\rho_i}(x)}|\nabla u_n(z)|^p\,dz}{\rho_i^{N}} \\
\le \frac{\rho_i^{p-1}}{\alpha} 
\left( \frac{\Es_n(u_n,B_{\rho_i}(x))}{\rho_i^{N-1}}
-\frac{\int_{B_{\rho_i}(x)}a_1(y)\,dy}{\rho_i^{N-1}}
\right).
\end{multline*}
Since $u_n$ is optimal for $u$ and by \eqref{ginftyux} we have that 
$$
\frac{\Es_n(u_n,B_{\rho_i}(x))}{\rho_i^{N-1}} \stackrel{n \to +\infty}\longrightarrow
\frac{\Es(u,B_{\rho_i}(x))}{\rho_i^{N-1}}
\stackrel{i \to +\infty}\longrightarrow
\omega_{N-1}g_\infty(x,u^-(x),u^+(x),\nu_{S(u)}(x))<+\infty.
$$
In view also of \eqref{a1rhoi}, we conclude that we can choose $n_i$ so that for $n \ge n_i$ 
\begin{equation*}
\label{gradstep3} 
\int_{B_1}|\nabla v^i_n(y)|^p\,dy \le 
C\rho_i^{p-1}
\end{equation*}
for some constant $C \ge 0$.
By Coarea formula for $BV$ functions (see \cite[Theorem 3.40]{AFP}) we get
$$
\int_{u^-(x)}^{u^+(x)} 
\hn \left( \partial^*E^i_n(t) \setminus S(v^i_n) \right) \,dt
\le
\int_{B_1}|\nabla v^i_n|\,dy \le \tilde C \rho_i^{1-\frac{1}{p}},
$$
for a suitable constant $\tilde C$, where 
\begin{equation*}
\label{eni}
E^i_n(t):=\{x \in B_1\,:\, x \text{ is a Lebesgue point for }v^i_n \text{ and }v^i_n(x)>t\}
\end{equation*}
and $\partial^*$ denotes the reduced boundary. By the Mean Value Theorem there exists $t^i_n \in [u^-(x),u^+(x)]$ such that
\begin{equation}
\label{meanvalue}
\hn \left( \partial^*E^i_n(t^i_n) \setminus S(v^i_n) \right)
\le \frac{\tilde C}{u^+(x)-u^-(x)} \rho_i^{1-\frac{1}{p}}.
\end{equation}
We now employ a construction similar to that employed by Francfort and Larsen in their Transfer of Jump Sets Theorem \cite[Theorem 2.3]{FL}. Since $x$ is a jump point for $u$ we have that for $i \to +\infty$ 
$$
u(x+\rho_i y) \to u_{0,u^-(x),u^+(x),\nu_{S(u)}(x)}
\qquad \text{strongly in }L^1(B_1).
$$
Then we have that for $n$ large
$$
| B_1^+ \bigtriangleup E^i_n(t^i_n) | \le e_i,
$$
where $B_1^+:=\{y \in B_1\,:\, y\cdot \nu_{S(u)}(x) \ge 0\}$,
$A \bigtriangleup B:=(A \setminus B) \cup (B \setminus A)$,
and $e_i \to 0$ for $i \to +\infty$. By Fubini's Theorem we have
$$
\int_0^{\sqrt{e_i}} \hn \left( (B_1^+ \setminus E^i_n(t^i_n)) \cap H(s) \right)\,ds \le
\int_{-\infty}^{+\infty} \hn \left( (B_1^+ \setminus E^i_n(t^i_n)) \cap H(s) \right)\,ds \le e_i,
$$
where $H(s):=\{y \in B_1\,:\, y\cdot \nu_{S(u)}(x)=s\}$, and by the Mean Value Theorem we get that there exists $0<s^{i,+}_n<\sqrt{e_i}$ such that setting $H^{i,+}_n:=H(s^{i,+}_n)$ we have
$$
\hn \left( (B_1^+ \setminus E^i_n(t^i_n)) \cap H^{i,+}_n \right) \le \sqrt{e_i}.
$$
Similarly we obtain $-\sqrt{e_i}<s^{i,-}_n<0$ such that setting $H^{i,-}_n:=H(s^{i,-}_n)$ we have
$$
\hn \left( (E^i_n(t^i_n) \setminus B_1^+) \cap H^{i,-}_n \right) \le \sqrt{e_i}.
$$
Let us write $y=(y',y_N)$, where $y_N$ is the coordinate along $\nu_{S(u)}(x)$ and $y'$ the coordinates in the hyperplain orthogonal to $\nu_{S(u)}(x)$. Let $l_i$
be such that for every $y \in B_1$
$$
|y_N| \ge 2\sqrt{e_i}
\Longrightarrow
|y'| \le 1-l_i.
$$
Let us set
$$
D^i_n:=(E^i_n(t^i_n)  \cup \{y \in B_1\,:\, y_N \ge s^{i,+}_n\}) \setminus 
\{y \in B_1\,:\, y_N \le s^{i,-}_n\}.
$$
We set
\begin{equation*}
\label{defwni}
w^i_n:=
\begin{cases}
1 & |y'| \ge 1-l_i, y_N \ge 0, \\
0 & |y'| \ge 1-l_i, y_N < 0, \\
1 & |y'| \le 1-l_i, y \in D^i_n, \\
0 & \text{otherwise.}
\end{cases}
\end{equation*}
Notice that $w^n_i$ is piecewise constant, with $w^i_n=u_{0,0,1,\nu_{S(u)}(x)}$ in a neighborhood of the boundary, and such that
\begin{equation*}
\label{ineqstep3bis}
\int_{B_1 \cap [S(v^i_n) \setminus K^i_n]}g_n(x+\rho_iy,\nu)\,d\hn(y)
\ge
\int_{B_1 \cap [S(w^i_n) \setminus K^i_n]}g_n(x+\rho_iy,\nu)\,d\hn(y)
-\tilde e_i, 
\end{equation*}
with $\tilde e_i \to 0$ for $i \to +\infty$.
\par
In view of \eqref{meanvalue} and of the assumption
\eqref{muxfinite} we have that $\hn(S(w^i_n)) \le C_i$ uniformly in $n$ for some finite constant $C_i$.  
By Ambrosio's Compactness Theorem (see for example \cite[Theorem 4.8]{AFP}) we get for $n \to +\infty$
$$
w^i_n \to w^i
\quad\quad
\text{strongly in }L^1(B_1),
$$
where $w^i$ is piecewise constant and $w^i=u_{0,0,1,\nu_{S(u)}(x)}$ in a neighborhood of the boundary, so that the claim is proved. 
\end{proof}

\begin{remark}
\label{nonintrem}
{\rm Theorem \ref{gammaresult} states that in the $\Gamma$-limit process there is no interaction between bulk and surface energies, since they are constructed looking at $\Gamma$-convergence problems in Sobolev space and in the space of piecewise constant functions respectively. As a consequence, considering 
bulk and surface energy densities of the form $c_1f_n$ and $c_2g_n$ with $c_1,c_2>0$, we get in the limit $c_1f$ and $c_2g$ as bulk and surface energy densities. We remark that a key assumption for non interaction is given by equi-boundness of $\hn(K_n)$: dropping this assumption, interaction can occur even in the case of constant densities, for example $f(\xi):=|\xi|^p$ and $g(x,\nu) \equiv 1$ (if we consider in $]0,1[$ the set $K_n:=\{\frac{i}{n}\,:\,i=1,\dots,n-1\}$, we get as $\Gamma$-limit the zero functional). As mentioned in the Introduction, non interaction between bulk and surface energies was noticed in the case of periodic homogenization (with $K_n=\emptyset$) by Braides, Defranceschi and Vitali in \cite{BDV}.
}
\end{remark}

In the rest of this section we employ Theorem \ref{gammaresult} to obtain 
a lower semicontinuity result for $SBV$ functions in the case of varying bulk and surface energies in the same spirit of Ambrosio's lower semicontinuity theorems \cite{A2}. 
\par
From Theorem \ref{gammaresult} we get that the following semicontinuity result holds.

\begin{proposition}
\label{lscprop2}
Let $(K_n)_{n \in \N}$ be a sequence of rectifiable sets in $\Om$ such that $\hn(K_n) \le C$ for all $n \in \N$.
Let us assume that for all $A \in \as(\Om)$ the functionals 
$\fs_n(\cdot,A)$ and $\gs^-_n(\cdot,A)$ defined in \eqref{bulkfuncn} and \eqref{surffuncnm} $\Gamma${-}converge in the strong topology of $L^1(\Om)$ to $\fs(\cdot,A)$ and $\gs^-(\cdot,A)$ respectively.
Let $(u_n)_{n \in \N}$ be a sequence in $SBV^p(\Om)$ such that 
$u_n \weak u$ weakly in $SBV^p(\Om)$. 
\par
Then for all $A \in \as(\Om)$ we have
\begin{equation}
\label{lsc1bis}
\int_A f(x,\nabla u(x))\,dx \le \liminf_{n \to +\infty}
\int_A f_n(x,\nabla u_n(x))\,dx,
\end{equation}
and
\begin{equation}
\label{lsc2bis}
\int_{S(u) \cap A} g^-(x,\nu)\,d\hn \le \liminf_{n \to +\infty} 
\int_{(S(u_n) \setminus K_n)\cap A} g_n(x,\nu)\,d\hn,
\end{equation}
where $f$ and $g^-$ are the densities of $\fs$ and $\gs^-$ respectively.
\par
In particular if $K_n=\emptyset$ we have
\begin{equation}
\label{lsc2}
\int_{S(u) \cap A} g(x,\nu)\,d\hn \le \liminf_{n \to +\infty} 
\int_{S(u_n) \cap A} g_n(x,\nu)\,d\hn,
\end{equation}
where $g$ is the density of $\gs$ defined in Proposition \ref{formulag2}.
\end{proposition}

\begin{proof}
By Theorem \ref{gammaresult}, we have that for all $h,k \in \N$ and for all $A \in \as(\Om)$ the functionals
$$
\Es_n^{h,k}(u,A):=
h\int_A f_n(x,\nabla u(x))\,dx+
k\int_{(S(u) \setminus K_n) \cap A} g_n(x,\nu)\,d\hn
$$
$\Gamma${-}converge in the strong topology of $L^1(\Om)$ to
$$
\Es^{h,k}(u,A):=
h\int_A f(x,\nabla u(x))\,dx+
k\int_{S(u) \cap A} g^-(x,\nu)\,d\hn.
$$
In particular by $\Gamma${-}liminf inequality we have
$$
\Es^{h,k}(u,A) \le \liminf_{n \to +\infty} \Es^{h,k}_n(u_n,A). 
$$
Then we get
\begin{multline*}
\int_A f(x,\nabla u(x))\,dx \le 
\liminf_{n \to +\infty} 
\int_A f_n(x,\nabla u_n(x))\,dx 
+\frac{k}{h}\int_{(S(u_n) \setminus K_n)\cap A}g_n(x,\nu)\,d\hn(x)  \\
\le \liminf_{n \to +\infty} 
\int_A f_n(x,\nabla u_n(x))\,dx 
+\frac{k}{h}C
\end{multline*}
for some constant $C$ independent of $h$ and $k$. 
Since $h,k$ are arbitrary we get that
\eqref{lsc1bis} holds. The proof of \eqref{lsc2bis} is analogous.
\end{proof}

\section{A new variational convergence for rectifiable sets}
\label{sigmaconv}
In this section we use the $\Gamma$-convergence results of Section \ref{gammasec} 
to introduce a variational notion of convergence for rectifiable sets which will be employed in the study of stability of unilateral minimality properties.
\par
Let $(K_n)_{n \in \N}$ be a sequence of rectifiable sets in $\Om$, and
let us assume following Ambrosio and Braides \cite[Theorem 3.2]{AB1} that the functionals
$\hs^-_n:P(\Om) \times \as(\Om) \to [0,+\infty)$ defined by
\begin{equation}
\label{hsnfunc}
\hs^-_n(u,A):=
\hn\left( (S(u) \setminus K_n) \cap A \right)  
\end{equation}
$\Gamma${-}converge with respect to the strong topology of $L^1(\Om)$ for every $A \in \as(\Om)$ to a functional $\hs^-(\cdot,A)$, which by the representation result of 
Bouchitt\'e, Fonseca, Leoni and Mascarenhas \cite[Theorem 3]{BFLM} is of the form
\begin{equation}
\label{hsfunc}
\hs^-(u,A):=
\int_{S(u) \cap A}h^-(x,\nu)\,d\hn(x)  
\end{equation}
for some function $h^-:\Om \times S^{N-1} \to [0,+\infty)$.

\begin{definition}[\bf $\sigma$-convergence of rectifiable sets]
\label{sigmadef}
Let $(K_n)_{n \in \N}$ be a sequence of rectifiable sets in $\Om$. We say that $K_n$  $\sigma$-converges in $\Om$ to $K$ if the functionals $(\hs^-_n)_{n \in \N}$ defined in \eqref{hsnfunc} $\Gamma$-converge in the strong topology of $L^1(\Om)$ to the functional $\hs^-$ defined in \eqref{hsfunc}, and $K$ is the (unique) rectifiable set in $\Om$ such that
\begin{equation}
\label{hzerok}
h^-(x,\nu_K(x))=0 \text{ for $\hn$-a.e. }x \in K,
\end{equation}
and such that for every rectifiable set $H \subseteq \Om$ we have
\begin{equation}
\label{hzerokprop}
h^-(x,\nu_H(x))=0 \text{ for $\hn$-a.e. }x \in H \Longrightarrow H \tsub K,
\end{equation}
where $H \tsub K$ means that $H \subseteq K$ up to a set of $\hn$-measure zero.
\end{definition}

\begin{remark}
\label{perturbationlem}
{\rm
From Definition \ref{sigmadef} it comes directly that $\sigma$-convergence of rectifiable sets is stable under infinitesimal perturbation in surface. More precisely,
let $(K_n)_{n \in \N}$ be a sequence of rectifiable sets in $\Om$ such that $K_n$ $\sigma$-converges in $\Om$ to $K$, and let $(\tilde K_n)_{n \in \N}$ be a sequence
of rectifiable sets in $\Om$ such that $\hn( \tilde K_n \bigtriangleup K_n) \to 0$, where
$\bigtriangleup$ denotes the symmetric difference of sets.
Then $\tilde K_n$ $\sigma$-converges in $\Om$ to $K$.
}
\end{remark}

Let us now come to the main properties of $\sigma$-convergence for rectifiable sets.
By compactness of $\Gamma$-convergence, we deduce the following compactness result for $\sigma$-convergence.

\begin{proposition}[\bf compactness]
\label{sigmacomp}
Let $(K_n)_{n \in \N}$ be a sequence of rectifiable sets in $\Om$ with \linebreak
\hbox{$\hn(K_n) \le C$.}
Then there exists a subsequence $(n_h)_{h \in \N}$ and a rectifiable set $K$ in $\Om$ such that $K_{n_h}$ $\sigma$-converges in $\Om$ to $K$. Moreover
\begin{equation}
\label{kliminf}
\hn(K) \le \liminf_{n \to +\infty} \hn(K_n).
\end{equation}
\end{proposition}

\begin{proof}
By Proposition \ref{formulag}, up to a subsequence we have that for all $A \in \as(\Om)$ the functionals $\hs^-_n(\cdot,A)$ defined in \eqref{hsnfunc} $\Gamma$-converge in the strong topology of $L^1(\Om)$ to a functional $\hs^-(\cdot,A)$ which can be represented through a density $h^-$ according to \eqref{hsfunc}.
\par
Let us consider the class
$$
\ks:=\{H \subseteq \Om: H \text{ is rectifiable} \text{ and } h^-(x,\nu_H(x))=0 \text{ for $\hn$-a.e. }x \in H\}.
$$
Notice that $\ks$ contains at least the empty set. Let us prove that for all $H \in \ks$ we have
\begin{equation}
\label{sigmahnh2}
\hn(H) \le L:=\liminf_{n \to +\infty}\hn(K_n).
\end{equation}
In fact let $H \in \ks$. Since $H=\cup_i H_i$ with $H_i$ compact and rectifiable with
$\hn(H_i)<+\infty$, it is not restrictive to consider $\hn(H)<+\infty$. Given $\eps>0$, by a covering argument we can find an open set $U$ and a piecewise constant function $v \in P(\Om)$ such that 
$$
\hn(H \setminus U)<\eps
\quad \text{ and }\quad
\hn \left( (S(v) \bigtriangleup H) \cap U \right) < \eps,
$$ 
where $\bigtriangleup$ denotes the symmetric difference of sets. Since $h^- \le 1$ we have
$$
\hs^-(v,U)=\int_{S(v) \cap U}h^-(x,\nu)\,d\hn(x)=
\int_{(S(v) \setminus H) \cap U}h^-(x,\nu)\,d\hn(x) < \eps.
$$
Let $(v_n)_{n \in \N}$ be a recovering sequence for $v$ with respect to $\hs^-(\cdot,U)$.
Then we have that 
$$
\limsup_{n \to +\infty} \hn\left( (S(v_n) \setminus K_n) \cap U \right)<\eps.
$$
By Ambrosio's Theorem we deduce that
\begin{multline*}
\hn(H) \le 
\hn(H \cap U)+\hn(H \setminus U) \le 
\hn(S(v) \cap U)+2\eps 
\\
\le \liminf_{n \to +\infty}\hn(S(v_n) \cap U)+2\eps 
\le \liminf_{n \to +\infty} \hn(K_n)+3\eps
= L+3\eps.
\end{multline*}
Since $\eps$ is arbitrary we get that \eqref{sigmahnh2} holds. 
\par
Let us now consider
$$
\tilde L:=\sup\{ \hn(H)\,:\, H \in \ks\}<+\infty,
$$
and let $(H_k)_{k \in \N}$ be a maximizing sequence for $\tilde L$.
%We can assume that $H_k \tsub H_m$ for $k \le m$. 
We set
$K:= \bigcup_{k=1}^\infty H_k$.
Clearly \eqref{kliminf} and \eqref{hzerok} hold. Moreover, since $\hn(K)=\tilde L$ we have that \eqref{hzerokprop} holds, and the proof is concluded. 
\end{proof}

\begin{remark}
\label{sigmarem}
{\rm
Let $\Om:=(-1,1) \times (-1,1)$ in $\R^2$, and let $(K_n)_{n \in \N}$ be a sequence of closed sets with $K_n \to K:=\{(-1,1)\} \times \{0\}$ in the Hausdorff metric and such that
$\hs^1 \res K_n \weakst a \hs^1 \res K$
\wlystar in the sense of measures. 
If $a<1$ by \eqref{kliminf} we deduce that $K_n$ $\sigma$-converges in $\Om$ to the empty set. We stress that the condition $a \ge 1$ is not enough to guarantee that $K$ is the $\sigma$-limit of $(K_n)_{n \in \N}$. In fact considering
$$
K_n:= \bigcup_{i=-n}^n \left\{ \frac{i}{n} \right\} \times \left[-\frac{1}{n},\frac{1}{n} \right]
$$
we have
$\hs^1 \res K_n \weakst 2 \hs^1 \res K$
\wlystar in the sense of measures. However also in this case we have that $K_n$ $\sigma$-converges in $\Om$ to the empty set. In fact let us consider $u \in P(\Om)$ such that
$u=1$ in $\Om^+:=(-1,1) \times (0,1)$ and $u=0$ in $\Om^-:=(-1,1) \times (-1,0)$, and let $u_n$ be a sequence in $P(\Om)$ such that $u_n \to u$ strongly in $L^1(\Om)$ and with $\hn(S(u_n)) \le C$. Let $(e_1,e_2)$ be the canonical base of $\R^2$.
By Ambrosio's theorem we get that
$$
\nu [u_n] \hs^1 \res S(u_n) \weakst
e_2 \hs^1 \res S(u)
$$
\wlystar in the sense of measures. Considering the vector field $\varphi e_2$ with $\varphi \in C^\infty_c(\Om)$
we get
\begin{equation*}
\int_{S(u_n) \setminus K_n} \varphi e_2 \cdot \nu [u_n] \,d\hs^1=
\int_{S(u_n)} \varphi e_2 \cdot \nu [u_n] \,d\hs^1 \to \int_K \varphi \,d\hs^1.
\end{equation*}
Since $\varphi$ is arbitrary, we deduce that $\liminf_{n \to +\infty} \hs^-_n(u_n)=\liminf_{n \to +\infty} \hs^1(S(u_n) \setminus K_n) \ge 1$. By $\Gamma$-liminf we conclude that $\hs^-(u)=1$ that is $h^-(x,e_2)=1$ for $\hs^1$-a.e.
$x \in K$. Since the $\sigma$-limit of $(K_n)_{n \in \N}$ can be only contained in $K$, we deduce that the $\sigma$-limit is the empty set.
}
\end{remark}

The following proposition, which comes immediately from the growth estimates on $g_n$, shows that the $\sigma$-limit is a natural limit candidate for a sequence of rectifiable sets in connection with unilateral minimality properties (see the Introduction).

\begin{proposition}
\label{sigmaprop1}
Let $(K_n)_{n \in \N}$ be a sequence of rectifiable sets in $\Om$ with 
$K_n$ $\sigma$-converging in $\Om$ to $K$.
Let $(g_n)_{n \in \N}$ be a sequence of Borel functions satisfying the growth estimates \eqref{surfenergy}, and let $g^-$ be the energy density of the $\Gamma$-limit in the strong topology of $L^1(\Om)$ of the functionals $(\gs^-_n)_{n \in \N}$ defined in \eqref{surffuncnm}. Then we have
\begin{equation*}
\label{gmzerok}
g^-(x,\nu_K(x))=0 \text{ for $\hn$-a.e. }x \in K,
\end{equation*}
and for every rectifiable set $H \subseteq \Om$ 
\begin{equation*}
\label{gmzerokprop}
g^-(x,\nu_H(x))=0 \text{ for $\hn$-a.e. }x \in H \Longrightarrow H \tsub K.
\end{equation*}
\end{proposition}

The following lower semicontinuity result for surface energies along sequences of rectifiable sets converging in the sense of $\sigma$-convergence will be employed in Section \ref{Composites}.

\begin{proposition}[\bf lower semicontinuity]
\label{lscsigmap}
Let $(K_n)_{n \in \N}$ be a sequence of rectifiable sets in $\Om$ 
such that $K_n$ $\sigma$-converges in $\Om$ to $K$. Let $(g_n)_{n \in \N}$ be a sequence of Borel functions satisfying the growth estimates \eqref{surfenergy}, and let $g$ be the associated function according to Proposition \ref{formulag2}. Then we have
$$
\int_K g(x,\nu)\,d\hn(x) \le \liminf_{n \to +\infty}
\int_{K_n} g_n(x,\nu)\,d\hn(x).
$$
\end{proposition}

\begin{proof}
Let $H \tsub K$ with $\hn(H)<+\infty$. Given $\eps>0$,  by a covering argument 
we can find an open set $U$ and a piecewise constant function $v \in P(\Om)$ such that 
$$
\hn(H \setminus U)<\eps
\quad\text{ and }\quad
\hn \left( (S(v) \bigtriangleup H) \cap U \right) < \eps,
$$ 
where $\bigtriangleup$ denotes the symmetric difference of sets. 
If $(v_n)_{n \in \N}$ is a recovering sequence for $v$ with respect to $\hs^-(\cdot,U)$
defined in \eqref{hsfunc} we have
$$
\limsup_{n \to +\infty} \hn\left( (S(v_n) \setminus K_n) \cap U \right)<\eps.
$$
We deduce by $\Gamma$-convergence that
\begin{multline*}
\int_{H} g(x,\nu)\,d\hn(x)=
\int_{H \cap U} g(x,\nu)\,d\hn(x)+
\int_{H \setminus U} g(x,\nu)\,d\hn(x) \\
\le \int_{S(v) \cap U} g(x,\nu)\,d\hn(x)+2\beta\eps 
\le \liminf_{n \to +\infty}\int_{S(v_n) \cap U} g_n(x,\nu)\,d\hn(x)+2\beta\eps \\
\le \liminf_{n \to +\infty}\int_{K_n} g_n(x,\nu)\,d\hn(x)+3\beta \eps.
\end{multline*}
Since $\eps$ is arbitrary we deduce
$$
\int_{H} g(x,\nu)\,d\hn(x) \le
\liminf_{n \to +\infty}\int_{K_n} g_n(x,\nu)\,d\hn(x),
$$
and since $H$ is arbitrary in $K$ the proof is concluded.
\end{proof}

The following proposition is essential in the study of stability of unilateral minimality properties.

\begin{proposition}
\label{sigmaprop2}
Let $(K_n)_{n \in \N}$ be a sequence of rectifiable sets in $\Om$ such that $K_n$ $\sigma$-converges in $\Om$ to $K$. Let $1<p<+\infty$, and let $(u_n)_{n \in \N}$ be a sequence in $SBV^p(\Om)$ with $u_n \weak u$ weakly in $SBV^p(\Om)$ and $\hn(S(u_n) \setminus K_n) \to 0$. Then $S(u) \tsub K$.
\end{proposition}

\begin{proof}
Let us consider $\tilde K_n:= S(u_n) \cap K_n$. By compactness, up to a further subsequence we have that $\tilde K_n$ $\sigma$-converges in $\Om$ to a rectifiable set $\tilde K \tsub K$.
Let $\tilde h^-$ be the density associated to $(\tilde K_n)_{n \in \N}$ according to Definition \ref{sigmadef}. By lower semicontinuity given by Proposition \ref{lscprop2} we have
$$
\int_{S(u)} \tilde h^-(x,\nu)\,d\hn(x) \le \liminf_{n \to +\infty} \hn \left( S(u_n) \setminus \tilde K_n \right)=0.
$$
We deduce that
$$
\tilde h^-(x,\nu_{S(u)}(x))=0 \text{ for $\hn$-a.e. }x \in S(u),
$$
so that by definition of $\sigma$-limit we deduce $S(u) \tsub \tilde K \tsub K$.
\end{proof}

The next corollary shows that our $\sigma$-limit always contains the $\sigma^p$-limit of introduced by Dal Maso, Francfort and Toader in \cite{DMFT} to study quasistatic crack growth in nonlinear elasticity. 
We recall that $K_n$ $\sigma^p$-converges in $\Om$ to $K$ if the following hold:
\begin{itemize}
\item[(1)] if $u_h \weak u$ weakly in $SBV^p(\Om)$ with $S(u_h) \subseteq K_{n_h}$, then $S(u) \subseteq K$;
\vskip4pt
\item[(2)] $K=S(u)$ and there exists  $u_n \weak u$ weakly in $SBV^p(\Om)$ with $S(u_n) \subseteq K_n$.
\end{itemize}

\begin{corollary}
\label{sigmasigmap}
Let $(K_n)_{n \in \N}$ be a sequence of rectifiable sets in $\Om$ such that
$K_n$ $\sigma$-converges in $\Om$ to $K$. Let $1<p<+\infty$, and let us assume that $K_n$ $\sigma^p$-converges in $\Om$ to some rectifiable set $\tilde K$. Then $\tilde K \tsub K$.
\end{corollary}

\begin{proof}
The proof readily follows from Proposition \ref{sigmaprop2} and point $(2)$ of the definition of $\sigma^p$-convergence.
\end{proof}

\begin{remark}
\label{kstrict}
{\rm
Notice that is general we can have that the $\sigma^p$-limit $\tilde K$ of $(K_n)_{n \in \N}$ is strictly contained in $K$. In fact we can consider $\Om:=(-1,1) \times (-1,1)$ in $\R^2$, and 
$$
K_n:=\{(-1,1) \setminus L_n\} \times \{0\}
$$
with $L_n \subseteq (-1,1)$ and $|L_n| \to 0$. In this case we get $K=(-1,1) \times \{0\}$, while if $L_n$ is chosen in such a way that its $c_p$-{\it capacity} is big enough (see the celebrated example of the Neumann sieve, we refer to \cite{M}) we get $\tilde K=\emptyset$.
\par
This example is based on the fact that the $\sigma^p$-limit is influenced by infinitesimal perturbations of the $K_n$ as pointed out in Remark \ref{perturbationlem}, while the set $K$ is not. 
}
\end{remark}

In Section \ref{stabilitysecbdry} and Section \ref{Composites}, we will need a definition of $\sigma$-convergence in the closed set $\Omb$.

\begin{definition}[\bf $\sigma$-convergence in $\Omb$]
\label{sigmabardef}
Let $(K_n)_{n \in \N}$ be a sequence of rectifiable sets in $\Omb$. We say that $K_n$  $\sigma$-converges in $\Omb$ to $K \tsub \Omb$ if $K_n$ $\sigma$-converges in $\Om'$ to $K$ for every open bounded set $\Om'$ such that $\Omb \subseteq \Om'$.
\end{definition}

Notice that to check the $\sigma$-convergence in $\Omb$ of rectifiable sets, it is enough check $\sigma$-convergence in $\Om'$ for just one $\Om'$ with $\Omb \subseteq \Om'$.

\section{Stability of unilateral minimality properties}
\label{stabilitysec}
In this section we apply the results of Section \ref{gammasec} and Section \ref{sigmaconv} to obtain the stability result of unilateral minimality properties under $\Gamma${-}convergence for bulk and surface energies.

\begin{definition}[\bf unilateral minimizers]
\label{defminunil}
Let $f: \Om \times \R^N \to [0,+\infty[$ be a Carath\'eodory function and let $g: \Om \times S^{N-1} \to [0,+\infty[$ be a Borel function satisfying the growth estimates \eqref{bulkenergy} and \eqref{surfenergy}. We say that the pair $(u,K)$ with $u \in SBV^p(\Om)$ and $K$ rectifiable set in $\Om$ is a {\it unilateral minimizer} with respect to $f$ and $g$ if
$S(u) \tsub K$, and
$$
\int_\Om f(x,\nabla u(x))\,dx
\le \int_\Om f(x,\nabla v(x))\,dx+\int_{H \setminus K}g(x,\nu),
$$
for all pairs $(v,H)$ with $v \in SBV^p(\Om)$, $H$ rectifiable set in $\Om$ such that $S(v) \tsub H$ and $K \tsub H$.
\end{definition}

As in the previous sections, let $f_n: \Om \times \R^N \to [0,+\infty[$ be a Carath\'eodory function and let $g_n: \Om \times S^{N-1} \to [0,+\infty[$ be a Borel function
satisfying the growth estimates \eqref{bulkenergy} and \eqref{surfenergy}.
\par
Let us assume that the functionals $(\fs_n(\cdot,A))_{n \in \N}$ and $(\gs_n(\cdot,A))_{n \in \N}$ 
defined in \eqref{bulkfuncn} and \eqref{surffuncn} $\Gamma${-}converge in the strong topology of $L^1(\Om)$ to $\fs(\cdot,A)$ and $\gs(\cdot,A)$ for every $A \in \as(\Om)$ respectively. Let $f$ be the density of $\fs$ according to Proposition
\ref{formulaf} and let $g$ be the density of $\gs$ according to Proposition 
\ref{formulag2}.
\par
The main result of the paper is the following stability result for unilateral minimality properties under $\sigma$-convergence of rectifiable sets (see Definition \ref{sigmadef}), and $\Gamma$-convergence of bulk and surface energies.

\begin{theorem}
\label{stabilitythm}
Let $(u_n)_{n \in \N}$ be a sequence in $SBV^p(\Om)$ with
$u_n \weak u$ weakly in $SBV^p(\Om)$, and let $(K_n)_{n \in \N}$ be a sequence of rectifiable sets in $\Om$ with $\hn(K_n) \le C$ and such that $K_n$ $\sigma$-converges in $\Om$ to $K$. Let us assume that the pair $(u_n,K_n)_{n \in \N}$ is a unilateral minimizer for $f_n$ and $g_n$. 
\par
Then $(u,K)$ is a unilateral minimizer for $f$ and $g$.
Moreover we have
\begin{equation}
\label{convfnf}
\lim_{n \to +\infty}\int_\Om f_n(x,\nabla u_n(x))\,dx=
\int_\Om f(x,\nabla u(x))\,dx.
\end{equation}
\end{theorem}

\begin{proof}
By Theorem \ref{gammaresult} we have that
the functionals
$$
\Es_n(u):=
\begin{cases}
\int_\Om f_n(x,\nabla u(x))\,dx+ 
\int_{S(u) \setminus K_n} g_n(x, \nu)\,d\hn(x)  
& u \in SBV^p(\Om), \\
+\infty & \text{otherwise}
\end{cases}
$$
$\Gamma${-}converge with respect to the strong topology of $L^1(\Om)$ to the functional
$$
\Es(u):=
\begin{cases}
\int_\Om f(x,\nabla u(x))\,dx+ 
\int_{S(u)} g^-(x, \nu)\,d\hn(x)  
& u \in SBV^p(\Om), \\
+\infty & \text{otherwise},
\end{cases}
$$
where $f$ and $g^-$ are defined in \eqref{reprf} and \eqref{reprgm} respectively, with $g^- \le g$.
\par
By Proposition \ref{sigmaprop2} we have $S(u) \tsub K$, so that $u$ is admissible for $K$, while by Proposition \ref{sigmaprop1} we have that
$$
g^-(x,\nu_K(x))=0 \text{ for $\hn$-a.e. }x \in K.
$$
Then the unilateral minimality of the pair $(u,K)$ easily follows. In fact, by $\Gamma$-convergence we have that $u$ is a minimizer for $\Es$ and $\Es_n(u_n) \to \Es(u)$. Then for all pairs $(v,H)$ with $S(v) \tsub H$ and $K \tsub H$ we have
\begin{multline*}
\int_\Om f(x,\nabla u(x))\,dx=\Es(u) \le \Es(v)=
\int_\Om f(x,\nabla v(x))\,dx+\int_{S(v)}g^-(x,\nu)\,d\hn \\
=\int_\Om f(x,\nabla v(x))\,dx+\int_{S(v) \setminus K}g^-(x,\nu) 
\le \int_\Om f(x,\nabla v(x))\,dx+\int_{H \setminus K}g(x,\nu),
\end{multline*}
so that the unilateral minimality property holds. The convergence of bulk energies
\eqref{convfnf} is given by the convergence $\Es_n(u_n) \to \Es(u)$.
\end{proof}

\begin{remark}[\bf stability under $\sigma^p$-convergence]
\label{stabilitysigmap}
{\rm
In the case of fixed bulk and surface energy densities $f$ and $g$, Dal Maso, Francfort and Toader \cite{DMFT} proved the stability of the unilateral minimality property under $\sigma^p$-convergence for the rectifiable sets $K_n$ (see Section \ref{sigmaconv} just before Corollary \ref{sigmasigmap} for the definition). This result readily follows by Theorem \ref{stabilitythm}. In fact
by Corollary \ref{sigmasigmap} we have that if $K_n$ $\sigma^p$-converges in $\Om$ to $\tilde K$, then $\tilde K$ is contained in the $\sigma$-limit of $(K_n)_{n \in \N}$. Since $S(u) \tsub \tilde K$, we get that the unilateral minimality of the pair $(u, \tilde K)$ is implied by the unilateral minimality of $(u,K)$.
}
\end{remark}

As mentioned in the Introduction, a method for proving stability of unilateral minimality properties nearer to the approach of \cite{DMFT} would be to prove a 
generalization of the Transfer of Jump Sets by Francfort and Larsen \cite[Theorem 2.1]{FL} to the case of varying energies. The following theorem based on the arguments of Section \ref{gammasec} provides such a generalization.

\begin{theorem}[\bf Transfer of Jump Sets]
\label{transferofjump}
Let $(K_n)_{n \in \N}$ be a sequence of rectifiable sets in $\Om$ with $\hn(K_n) \le C$ and $K_n$ $\sigma$-converging in $\Om$ to $K$.
For every $v \in SBV^p(\Om)$ there exists $(v_n)_{n \in \N}$ sequence in $SBV^p(\Om)$ with $v_n \weak v$ weakly in $SBV^p(\Om)$ and
such that
$$
\lim_{n \to +\infty} \int_\Om f_n(x,\nabla v_n(x))\,dx =
\int_\Om f(x,\nabla v(x))\,dx
$$
and
$$
\limsup_{n \to +\infty} \int_{S(v_n) \setminus K_n} g_n(x,\nu)\,d\hn(x) \le
\int_{S(v) \setminus K} g(x,\nu)\,d\hn(x).
$$
\end{theorem}

\begin{proof}
Let $(v_n)_{n \in \N}$ be a recovering sequence for $v$ with respect to $(\Es_n)_{n \in \N}$ defined in \eqref{funcn2}. By growth estimates on $f_n$ and $g_n$, and since $\hn(K_n) \le C$, we get $v_n \weak v$ weakly in $SBV^p(\Om)$. Since no interaction between bulk and surface energies occurs in view of Theorem \ref{gammaresult}, we get that
$$
\lim_{n \to +\infty} \int_\Om f_n(x,\nabla v_n(x))\,dx =
\int_\Om f(x,\nabla v(x))\,dx
$$
and
$$
\lim_{n \to +\infty}\int_{S(v_n) \setminus K_n} g_n(x,\nu)\,d\hn
=\int_{S(v)} g^-(x,\nu)\,d\hn \le \int_{S(v) \setminus K} g(x,\nu)\,d\hn
$$
because $g^-=0$ on $K$, and $g^- \le g$. 
\end{proof}

\section{Stability of unilateral minimality properties with boundary conditions}
\label{stabilitysecbdry}
In view of the application of Section \ref{Composites}, we need a stability result for unilateral minimality properties with boundary conditions.
\par
Let $\partial_D \Om \subseteq \partial \Om$.
In order to take into account a boundary datum on $\partial_D \Om$, we will use
the following notation:  if $u,\psi \in SBV(\Om)$ we set
\begin{equation}
\label{sgjump}
\Sg{\psi}{u}:=S(u) \cup 
\{x \in \partial_D \Om\,:\,u(x) \not= \psi(x)\},
\end{equation}
where the inequality on $\partial_D \Om$ is intended in
the sense of traces. 
\par
In order to set the problem, let $f_n: \Om \times \R^N \to [0,+\infty[$ be a Carath\'eodory function satisfying the growth estimate \eqref{bulkenergy}, and let $g_n: \Omb \times S^{N-1} \to [0,+\infty[$ be a Borel function satisfying the growth estimate \eqref{surfenergy}.
We consider unilateral minimality properties of the form
\begin{equation*}
\label{minbdry}
\int_\Om f_n(x, \nabla u_n)\,dx \le
\int_\Om f_n(x, \nabla v)\,dx+
\int_{H \setminus K_n}  g_n(x,\nu)\,d\hn(x) 
\end{equation*}
for every $v \in SBV^p(\Om)$ and for every rectifiable set $H$ in $\Omb$ such that $\Sg{\psi_n}{v} \tsub H$. Here $(K_n)_{n \in \N}$ is a sequence of rectifiable sets in $\Omb$ with $\hn(K_n) \le C$, $(u_n)_{n \in \N}$ is a sequence in $SBV^p(\Om)$ with $\Sg{\psi_n}{u_n} \tsub K_n$,
$\psi_n \in W^{1,p}(\Om)$ with $\psi_n \to \psi$ strongly in $W^{1,p}(\Om)$, and $\Sg{\psi_n}{\cdot}$ is defined in \eqref{sgjump}.
\par
In order to treat $\Sg{\psi_n}{\cdot}$ as an internal jump and in order to recover the surface energy on $\partial_D \Om$ for the minimality property in the limit, let us consider an open bounded set $\Om'$ such that $\Omb \subset \Om'$ and let us consider $g'_n:\Om' \times S^{N-1} \to [0,+\infty[$ such that
\begin{equation*}
\label{defgnprimo}
g'_n(x,\nu):=
\begin{cases}
g_n(x,\nu) &\text{if }x \in \Omb, \\
\beta+1 & \text{otherwise}.
\end{cases}
\end{equation*}
Let us consider the functionals $\gs'_n:P(\Om') \times \as(\Om') \to [0,+\infty]$
defined by
$$
\gs'_n(v,A):=
\int_{S(v) \cap A} g'_n(x,\nu)\,d\hn(x)
$$
and let $\gs':P(\Om') \times \as(\Om') \to [0,+\infty]$ be their $\Gamma${-}limit in the strong topology of $L^1(\Om')$, which according to Proposition \ref{formulag} is
of the form
\begin{equation}
\label{reprgprimo}
\gs'(v,A):=
\int_{S(v) \cap A} g'(x,\nu)\,d\hn(x).
\end{equation}
We clearly have $g'(x,\nu)=g(x,\nu)$ for $x \in \Om$, where $g$ is the surface energy density defined in \eqref{reprg2}, while it turns out that (see Remark \ref{gprimorem}) the surface energy given by the restriction of $g'$ to $\partial \Om \times S^{N-1}$ is completely determined by the functions $g_n$. 
\par
Let us set
\begin{equation*}
\label{deffnprimo2}
f'_n(x, \xi):=
\begin{cases}
f_n(x,\xi)      & \text{if }x \in \Om, \\
\alpha |\xi|^p      & \text{otherwise,}
\end{cases}
\end{equation*}
and let $f'$ be the energy density of the $\Gamma$-limit of the functionals on $W^{1,p}(\Om')$ associated to $f'_n$ according to Proposition \ref{formulaf}. We easily have that
\begin{equation*}
\label{deffprimo}
f'(x, \xi):=
\begin{cases}
f(x,\xi)      & \text{if }x \in \Om, \\
\alpha |\xi|^p      & \text{otherwise.}
\end{cases}
\end{equation*}
Since $\Om$ is Lipschitz, we can assume using an extension operator that $\psi_n,\psi \in W^{1,p}(\R^N)$ and $\psi_n \to \psi$ strongly in $W^{1,p}(\R^N)$.
\par
Before stating our stability result, we need the following $\Gamma$-convergence result,
which is a version of Theorem \ref{gammaresult} taking into account boundary data.

\begin{lemma}
\label{gammaconvlem}
Let $(K_n)_{n \in \N}$ be a sequence of rectifiable sets in $\Omb$ such that $\hn(K_n) \le C$. Let us assume that the functionals
\begin{equation*}
\label{esntilde2}
\Es'_n(v):=
\begin{cases}
\int_{\Om'}f'_n(x,\nabla v(x))\,dx+\int_{S(v) \setminus K_n}g'_n(x,\nu)\,d\hn(x)      
& \text{if }v \in SBV^p(\Om'),\\
+\infty      & \text{otherwise}
\end{cases}
\end{equation*}
$\Gamma${-}converge in the strong topology of $L^1(\Om')$ according to
Theorem \ref{gammaresult} to
\begin{equation*}
\label{estilde2}
\Es'(v):=
\begin{cases}
\int_{\Om'}f'(x,\nabla v(x))\,dx+\int_{S(v)}g'^{\,-}(x,\nu)\,d\hn(x)      
& \text{if }v \in SBV^p(\Om'), \\
+\infty      & \text{otherwise.}
\end{cases}
\end{equation*}
Then we have that the functionals
\begin{equation*}
\label{esntilde3}
\tilde \Es'_n(v):=
\begin{cases}
\Es'_n(v)      & \text{if }v=\psi_n \text{ on }\Om' \setminus \Omb, \\
+\infty      & \text{otherwise}
\end{cases}
\end{equation*}
$\Gamma${-}converge in the strong topology of $L^1(\Om')$ to
\begin{equation*}
\label{estilde3}
\tilde \Es'(v):=
\begin{cases}
\Es'(v)      & \text{if }v=\psi \text{ on }\Om' \setminus \Omb, \\
+\infty      & \text{otherwise.}
\end{cases}
\end{equation*}
\end{lemma}

\begin{proof}
\par
Let $v \in SBV^p(\Om')$ with $v=\psi$ on $\Om' \setminus \Omb$, and let $(v_n)_{n \in \N}$ be a recovering sequence for $v$ with respect to the functionals $\Es'_n$. 
We have that 
\begin{equation}
\label{vntopsi}
\nabla v_n \to \nabla \psi 
\qquad
\text{strongly in }L^p(\Om' \setminus \Omb;\R^N),
\end{equation}
and 
\begin{equation}
\label{svntozero}
\hn(S(v_n) \cap (\Om' \setminus \Omb)) \to 0.
\end{equation}
In fact we have that for all $U \in \as(\Om')$ such that $\overline U \subseteq \Om' \setminus \Omb$ and $\Es'(v,\partial U)=0$
\begin{equation*}
\label{vntopsiu}
\nabla v_n \to \nabla \psi 
\qquad
\text{strongly in }L^p(U;\R^N),
\end{equation*}
and
\begin{equation}
\label{svnuzero}
\hn(S(v_n) \cap U) \to 0.
\end{equation}
Let $\eps>0$ and let us consider an open set $V \in \as(\Om')$ such that $\partial \Om \subseteq V$, $\Es'(v,\partial V)=0$, $\int_{V \cap \Om}|a_1|\,dx<\eps$ ($a_1$ is defined in \eqref{bulkenergy}),
\begin{equation}
\label{bulku}
\int_V f'(x,\nabla v(x))\,dx < \eps
\quad\text{ and }\quad
\int_V f'(x,\nabla \psi(x))\,dx < \eps.
\end{equation}
Then for $n$ large (no interaction between bulk and surface part occurs) we have
\begin{equation}
\label{bulknu}
\int_V f'_n(x,\nabla v_n(x))\,dx < \eps.
\end{equation}
Notice that
\begin{multline*}
\int_{\Om' \setminus \Omb}|\nabla v_n- \nabla \psi|^p\,dx
= \int_{\Om' \setminus (\Om \cup V)}|\nabla v_n- \nabla \psi|^p\,dx +
\int_{V \setminus \Omb}|\nabla v_n- \nabla \psi|^p\,dx 
\\ \le
\int_{\Om' \setminus (\Om \cup V)}|\nabla v_n- \nabla \psi|^p\,dx +
\frac{2^{p-1}}{\alpha}\int_{V}f'_n(x,\nabla v_n(x))+f'(x,\nabla \psi(x))\,dx+
\frac{2^{p-1}}{\alpha}\int_{V \cap \Om}2|a_1|\,dx.
\end{multline*}
Since $\nabla v_n \to \nabla \psi$ strongly in $L^p(\Om' \setminus (\Om \cup V);\R^N)$, because of \eqref{bulku} and \eqref{bulknu}, and since $\eps$ is arbitrary, we get that \eqref{vntopsi} holds.
\par
Let us come to \eqref{svntozero}. Up to a subsequence we have
$$
\mu_n:=\hn \res (S(v_n) \cap (\Om' \setminus \Omb))) \weakst \mu
\qquad
\text{\wlystar in }\ms_b(\Om').
$$
In view of \eqref{svnuzero}, in order to prove \eqref{svntozero} it is sufficient to show that $\mu(\partial \Om)=0$. Let us assume by contradiction that $\mu(\partial \Om) \not= 0$: then there exists a cube $Q_\rho$ of center $x \in \partial \Om$ and edge $2\rho$
such that $\Es'(v,\partial Q_\rho)=0$ and
\begin{equation}
\label{densitymun}
\mu(Q_\rho)>\sigma>0.
\end{equation}
Up to a translation we may assume that $x=0$, and moreover we can assume that 
$$
\Om \cap Q_\rho=\{(x',y): x' \in (-\rho,\rho), y \in(-\rho,h(x'))\}, 
$$
where $(x',y)$ is a suitable orthogonal coordinate system and $h$ is a Lipschitz function.
Let $\eta>0$ be such that setting
$$
V_\eta:=\{(x',y):  x' \in (-\rho,\rho),\,\,  y \in (h(x')-\eta, h(x')+\eta) \}
$$
we have $V_\eta \subseteq Q_\rho$, and $\Es'(v,\partial V_\eta)=0$.
Let us set 
$$
V^-_\eta:=\{(x',y) \in V_\eta: y<h(x')\}
\quad
\text{ and }
\quad
V^+_\eta:=\{(x',y) \in V_\eta: y>h(x')\}.
$$
By \eqref{densitymun} we have that for $n$ large
\begin{equation}
\label{densitymun2}
\hn(S(v_n) \cap V^+_\eta)>\sigma.
\end{equation}
Let $\hat v$ be the function defined on $V_\eta$ obtained reflecting $v_{|V^+_\eta}$ to $V^-_\eta$: more precisely let us set
$$
\hat v=
\begin{cases}
v(x',y) & \text{if }(x',y) \in V^+_\eta, \\
v(x',2h(x')-y)   & \text{if }(x',y) \in V^-_\eta.
\end{cases}
$$
We clearly have $v \in W^{1,p}(V_\eta)$. Let $\hat v_n$ be obtained in the same way from $(v_n)_{|V^+_\eta}$. Let us consider 
$$
w_n:=v_n+\hat v-\hat v_n.
$$
We have $w_n \weak v$ weakly in $SBV^p(V_\eta)$ so that by lower semicontinuity given by Proposition \ref{lscprop2} we get
\begin{equation}
\label{veta1}
\int_{S(v) \cap V_\eta}g'^{\,-}(x,\nu)\,d\hn(x) \le
\liminf_{n \to +\infty} 
\int_{(S(w_n) \setminus K_n)\cap V_\eta}g'_n(x,\nu)\,d\hn(x).
\end{equation}
On the other hand, since $\Es'(v, \partial V_\eta)=0$, we have that $v_n$ is a recovering sequence for $v$ in $V_\eta$. In particular we get that
\begin{equation}
\label{veta2}
\int_{S(v) \cap V_\eta}g'^{\,-}(x,\nu)\,d\hn(x)=
\lim_{n \to +\infty} 
\int_{(S(v_n) \setminus K_n)\cap V_\eta}g'_n(x,\nu)\,d\hn(x).
\end{equation}
Formulas \eqref{veta1} and \eqref{veta2} give a contradiction because for $n$ large by \eqref{densitymun2} and since $K_n \tsub \Omb$ and $S(w_n) \tsub \Omb \cap Q_\rho$ (recall that $g'_n(x,\nu)=\beta+1$ for $x \in \Om' \setminus \Omb$)
$$
\int_{(S(v_n) \setminus K_n) \cap V_\eta}g_n(x,\nu)\,d\hn(x)
-\int_{(S(w_n) \setminus K_n) \cap V_\eta}g_n(x,\nu)\,d\hn(x)>
\sigma.
$$
We conclude that \eqref{svntozero} holds.
\par
We are now in a position to prove the $\Gamma$-limsup inequality for $\tilde \Es'_n$ and $\tilde \Es'$ (the $\Gamma$-liminf is immediate from the $\Gamma$-convergence of $\Es'_n$ to $\Es'$ and the fact that the constraint is closed under the strong topology of $L^1(\Om)$).
Let $\eps>0$, and let $U \in \as(\Om')$ be such that $\partial \Om \subseteq U$,
$\Es'(v,\partial U)=0$, and
\begin{equation}
\label{fpiccolou}
\int_U f(x,\nabla v)\,dx <\eps.
\end{equation}
In view of \eqref{vntopsi} and \eqref{svntozero} we can find $\varphi_n \in SBV^p(\Om')$ such that $\varphi_n=\psi_n-v_n$ on
$\Om' \setminus \Omb$, $\varphi_n=0$ on $\Om \setminus U$ and 
\begin{align*}
&\varphi_n \to 0 \qquad\text{strongly in }L^1(\Om'), \\
&\nabla \varphi_n \to 0 \qquad\text{strongly in }L^p(\Om';\R^N), \\
&\hn(S(\varphi_n)) \to 0.
\end{align*}
%$\varphi_n$ can be constructed by reflection of $(\psi_n-v_n)_{|\Om' \setminus \Omb}$ 
%near $\partial \Om$. 
Let us consider
$$
\tilde v_n:=v_n+\varphi_n.
$$
We have $\tilde v_n=\psi_n$ on $\Om' \setminus \Omb$. Moreover
$$
\limsup_{n \to +\infty} \int_{S(\tilde v_n) \setminus K_n}g'_n(x,\nu)\,d\hn =
\limsup_{n \to +\infty} \int_{S(v_n) \setminus K_n}g'_n(x,\nu)\,d\hn,
$$
and using the growth estimate on $f'_n$
\begin{multline*}
\limsup_{n \to +\infty} 
\left| \int_{\Om'} f'_n(x,\nabla \tilde v_n(x))\,dx-
\int_{\Om'} f'_n(x,\nabla v_n(x))\,dx \right| \\
\le \limsup_{n \to +\infty}
\int_{U \cap \Om} f_n(x,\nabla \tilde v_n(x))+f_n(x,\nabla v_n(x))\,dx \\
\le
\limsup_{n \to +\infty}
\int_U a_2(x)\,dx+ \left(\frac{2^{p-1}}{\alpha}+1 \right) \int_U f_n(x,\nabla v_n(x))\,dx\\
+
\frac{2^{p-1}}{\alpha}\int_U|a_1|\,dx+2^{p-1}
\int_U |\nabla \varphi_n|^p\,dx.
\end{multline*}
By \eqref{fpiccolou} we get
$$
\limsup_{n \to +\infty} \int_U f_n(x,\nabla v_n(x))\,dx <\eps.
$$
Then we conclude
$$
\limsup_{n \to +\infty} 
\left| \int_{\Om'} f'_n(x,\nabla \tilde v_n(x))\,dx-
\int_{\Om'} f'_n(x,\nabla v_n(x))\,dx \right| \le e(\eps),
$$
with $e(\eps) \to 0$ as $\eps \to 0$. We deduce that
$$
\limsup_{n \to +\infty} \tilde \Es'(\tilde v_n) \le \tilde \Es'(v)+e(\eps),
$$
with $e(\eps) \to 0$ as $\eps \to 0$. Since $\eps$ is arbitrary, using a diagonal argument
we have that the $\Gamma$-limsup inequality is proved.
\end{proof}

\begin{remark}
\label{gprimorem}
{\rm
In view of Lemma \ref{gammaconvlem} we can prove that the surface energy determined by the restriction of $g'$ to $\partial \Om$ is actually independent of the choice of $\Om'$ and of the constant value $c'$ of $g'_n$ on $\Om' \setminus \Omb$ provided that $c'>\beta$. In fact $g'$ is the density of the surface energy of the $\Gamma$-limit in the strong topology of $L^1(\Om)$ of the functionals on $SBV^p(\Om')$ defined as
$$
\hat \Es'_n(v):=\int_{\Om'} f'_n(x,\nabla v(x))\,dx+\int_{S(v)}g'_n(x,\nu)\,d\hn(x).
$$
Following the proof of Lemma \ref{gammaconvlem} (for the functionals $\Es'_n$ with $K_n=\emptyset$), if $v=\psi$ outside $\Omb$, we can find $(v_n)_{n \in \N}$  recovering sequence for $v$ with respect to $(\hat \Es'_n, \Om', c')$ such that $v_n=\psi_n$ outside $\Omb$. Then if $\Om''$ is an open set such that $\Omb \subseteq \Om''$ we have that $(v_n)_{|\Om' \cap \Om''}$ is a recovering sequence also for $(\hat \Es'_n, \Om'' \cap \Om', c'')$, and we have
$$
\int_{S(v)}g'(x,\nu)\,d\hn=
\lim_{n \to +\infty}
\int_{S(v_n)}g_n(x,\nu)\,d\hn.
$$
We deduce that the surface energy given by the restriction of $g'$ to $\Omb \times S^{N-1}$ is determined only by the $g_n:\Omb \times S^{N-1} \to [0,+\infty]$.
}
\end{remark}

The stability result for unilateral minimality properties with boundary conditions under $\sigma$-convergence in $\Omb$ for rectifiable sets (see Definition \ref{sigmabardef}) and $\Gamma$-convergence of bulk and surface energies is the following.

\begin{theorem}
\label{stabilitybdrythm}
Let $\psi_n \in W^{1,p}(\Om)$ with $\psi_n \to \psi$ strongly in $W^{1,p}(\Om)$. Let $(u_n)_{n \in \N}$ be a sequence in $SBV^p(\Om)$ with
$u_n \weak u$ weakly in $SBV^p(\Om)$, and let $(K_n)_{n \in \N}$ be a sequence of rectifiable sets in $\Omb$ with $\hn(K_n) \le C$, such that $K_n$ $\sigma$-converges in $\Omb$ to $K$, and $\Sg{\psi_n}{u_n} \tsub K_n$. 
\par
Let us assume that the pair $(u_n,K_n)$ satisfies the unilateral minimality property \eqref{minbdry} with respect to $f_n$, $g_n$ and $\psi_n$. Then $(u,K)$ satisfies the unilateral minimality property with respect to $f$, $g$ and $\psi$, where $f$ is defined in \eqref{reprf} and $g$ is the restriction of $g'$ defined in \eqref{reprgprimo} to $\Omb \times S^{N-1}$. 
Moreover we have
\begin{equation*}
\label{convfnfbdry}
\lim_{n \to +\infty}\int_\Om f_n(x,\nabla u_n(x))\,dx=
\int_\Om f(x,\nabla u(x))\,dx.
\end{equation*}
\end{theorem}

\begin{proof}
Since the boundary datum $\psi_n$ is imposed just on $\partial_D \Om$, 
we can consider $\partial_N \Om:=\partial \Om \setminus \partial_D \Om$ as part of the cracks, that is we can replace in the unilateral minimality properties $K_n$ with $K'_n:=K_n \cup \partial_N \Om$. 
\par
It is easy to prove that $K'_n$ $\sigma$-converges in $\Omb$ to $K \cup \partial_N \Om$. Then the proof follows that of Theorem \ref{stabilitythm} employing the functionals $(\tilde \Es'_n)_{n \in \N}$ defined in Lemma \ref{gammaconvlem} with $K'_n$ in place of $K_n$. 
\end{proof}

\section{Quasistatic evolution of cracks in composite materials}
\label{Composites}
The aim of this section is to apply the stability results of Section \ref{stabilitysecbdry}
to the study the asymptotic behavior of crack evolutions relative to varying bulk and surface energy densities $f_n$ and $g_n$. As mentioned in the Introduction, this problem is inspired by the problem of crack propagation in composite materials. 
We restrict our analysis to the case of antiplanar shear, where the elastic body is an infinite cylinder. 
\par
Let us recall the result of Dal Maso, Francfort and Toader \cite{DMFT} about quasistatic crack evolution in nonlinear elasticity: it is a very general existence and approximation result concerning a variational theory crack propagation inspired by the variational model introduced by Francfort and Marigo in \cite{FM}. As already said, we consider the antiplanar case and for simplicity we neglect body and traction forces, and so we adapt the mathematical tools employed in \cite{DMFT} to this scalar setting.
\par
As in the previous sections, let $\Om \subset \R^N$ (which, for $N=2$ represents a section of the cylindrical hyperelastic body) be an open bounded set with Lipschitz boundary. 
The family of admissible cracks is the class of rectifiable subsets of $\Omb$, while the class of admissible displacements is given by the functional space $SBV^p(\Om)$, where $1<p< +\infty$. Let $\partial_D \Om$ be a subset of $\partial \Om$. Given $\psi \in W^{1,p}(\Om)$, we say that the displacement $u$ is admissible for the fracture $K$ and the boundary datum $\psi$ and we write $u \in AD(\psi,K)$ if
$S(v) \tsub K$ and $v=\psi$ on $\partial_D \Om \setminus K$.
This can be summarized by the notation $\Sg{\psi}{u} \tsub K$, where $\Sg{\psi}{\cdot}$ is defined in \eqref{sgjump}.
\par  
Let $f(x,\xi): \Om \times \R^N \to [0,+\infty[$ be a Carath\'eodory function which is convex and $C^1$ in $\xi$ for a.e. $x \in \Om$, and satisfies the growth estimate
\begin{equation}
\label{cgrowthb}
a_1(x)+\alpha |\xi|^p \le f(x,\xi) \le a_2(x)+\beta |\xi|^p,
\end{equation}
where $a_1,a_2 \in L^1(\Om)$ and $\alpha, \beta>0$. Let moreover $g: \Omb \times S^{N-1} \to [0,+\infty[$ be a Borel function such that
\begin{equation}
\label{cgrowths}
\alpha \le g(x,\nu) \le \beta.
\end{equation}
The total energy of a configuration $(u,K)$ is given by
$$
\Es(u,K):=\int_\Om f(x, \nabla u(x))\,dx+\int_K g(x, \nu) d\hn(x).
$$
We will usually refer to the first term as {\it bulk energy} of $u$ and we write
\begin{equation*}
\label{cbulk}
\Eb(u):= \int_\Om f(x, \nabla u(x))\,dx, 
\end{equation*}
while we will refer to the second term as {\it surface energy} of $K$ and we write
\begin{equation*}
\label{csurf}
\Esup(K):= \int_K g(x, \nu) d\hn(x).
\end{equation*}
\par
Let us consider now a time dependent boundary datum $\psi \in W^{1,1}\big([0,T];W^{1,p}(\Om)\big)$ (i.e. the function $t \to \psi(t)$ is absolutely continuous from $[0,T]$ to the Banach space $W^{1,p}(\Om)$, with summable time derivative, see for instance \cite{Br1}), such that
for all $t \in [0,T]$
\begin{equation}
\label{linftypsi}
\|\psi(t)\|_{L^\infty(\Om)} \le C.
\end{equation}
In \cite{DMFT} Dal Maso, Francfort and Toader proved the esistence of an {\it irreversible quasistatic crack evolution} in $\Om$ relative to the boundary displacement $\psi$, i.e. the existence of a map $t \to (u(t),K(t))$ where $u(t) \in AD(\psi(t),K(t))$, $\|u(t)\|_{L^\infty(\Om)} \le \|\psi(t)\|_\infty$ and such that the following three properties hold:
\begin{itemize}
\vskip4pt
\item[(1)] {\it irreversibility:}  $K(t_1) \tsub K(t_2)$ for all $0 \le t_1\le t_2\le T$; 
\vskip4pt
\item[(2)] {\it static equilibrium:}
$\Es(u(0),K(0)) \le \Es(v,K)$ for all $(v,K)$ such that $v \in AD(\psi(0),K)$, and
$$
\Es(u(t),K(t)) \le \Es(v,K)
\qquad
\text{ for all }K(t) \tsub K, \,v \in AD(\psi(t),K);
$$
\vskip4pt
\item[(3)]  {\it energy balance:} the function $t \to \Es(u(t),K(t))$  
is absolutely continuous and
$$
\frac{d}{dt} \Es(u(t),K(t))=  
\int_\Om \nabla_{\xi}f(x,\nabla u(t)) \nabla \dot{\psi}(t)\,dx,
$$
where $\dot \psi$ denotes the time derivative of $t \to \psi(t)$.
\end{itemize}
\par
For every $n \in \N$ let us consider admissible bulk and surface energy densities
$f_n: \Om \times \R^N \to \R$ and $g_n: \Om \times S^{N-1} \to [0,+\infty[$ 
for the model of \cite{DMFT} satisfying the growth estimates
\eqref{cgrowthb} and \eqref{cgrowths} uniformly in $n$.
Let us moreover assume that $f_n$ is such that for a.e. $x \in \Om$ and for all $M \ge 0$
\begin{equation}
\label{fnhyp}
|\nabla_\xi f_n(x,\xi^1_n)- \nabla_\xi f_n(x,\xi^2_n)| \to 0
\end{equation}
for all $\xi^1_n,\xi^2_n$ such that $|\xi|^1_n \le M$
$|\xi|^2_n \le M$ and $|\xi^1_n-\xi^2_n| \to 0$. We denote by $\Es_n$, $\Eb_n$ and $\Esup_n$ the total, bulk and surface energies associated to $f_n$ and $g_n$.
\par
Let $f$ and $g$ be the effective energy densities associated to $f_n$ and $g_n$ in the sense of Theorem \ref{stabilitybdrythm}, i.e. let $f$ be given by Proposition \ref{formulaf} and let $g$ be the restriction to $\Omb \times S^{N-1}$ of the function $g'$ defined in \eqref{reprgprimo}. Notice that by Theorem \ref{regolare} we have that the function $f(x,\cdot)$ is $C^1$: as it is also convex in $\xi$ and satisfies the growth estimate \eqref{cgrowthb}, we have that
$f$ and $g$ are admissible bulk and surface energy densities for the model of \cite{DMFT}.
\par
Let $t \to \psi_n(t)$ be a sequence of admissible time dependent boundary displacements satisfying \eqref{linftypsi} and such that 
$$
\psi_n \to \psi
\qquad
\text{strongly in }W^{1,1} \Big([0,T],W^{1,p}(\Om)\Big).
$$
Let $t \to (u_n(t),K_n(t))$ be a quasistatic evolution for the boundary datum $\psi_n$
relative to the energy densities $f_n$ and $g_n$ according to \cite{DMFT}.
The main result of this section is the following Theorem which asserts that the $\sigma$-limit in $\Omb$ of $K_n(t)$ (see Definition \ref{sigmabardef}) still determines a quasistatic crack growth with respect to the energy densities $f$ and $g$.

\begin{theorem}
\label{compevol}
There exists a quasistatic crack growth $t \to (u(t),K(t))$ relative to the energy densities $f$ and $g$ and the boundary datum $\psi$ such that up to a subsequence (not rabelled) the following hold:
\begin{itemize}
\item[(1)] for all $t \in [0,T]$
$$
K_n(t) \text{ $\sigma$-converges in $\Omb$ to }K(t),
$$
and there exists a further subsequence $n_k$ (depending possibly on $t$) such that
$$
u_{n_k}(t) \weak u(t)
\qquad
\text{weakly in }SBV^p(\Om);
$$
\vskip4pt
\item[(2)] for every $t \in [0,1]$ we have convergence of total energies
$$
\Es_n(u_n(t),K_n(t)) \to \Es (u(t), K(t)),
$$
and in particular separate convergence for bulk and surface energies, i.e.
$$
\Eb_n(u_n(t)) \to \Eb(u(t)) 
\quad\text{ and }\quad
\Esup(K_n(t)) \to \Esup(K(t)).
$$
\end{itemize}
\end{theorem}

\begin{proof}
Notice that by the energy balance condition for $t \to (u_n(t),K_n(t))$ and by growth estimates on $f_n$ and $g_n$ we have that there exists a constant $C$ such that for all $t \in [0,T]$ and for all $n \in \N$
\begin{equation}
\label{cboundn}
\|\nabla u_n(t)\|^p+\hn(K_n(t))+\|u_n(t)\|_{L^\infty(\Om)} \le C.
\end{equation}
We divide the proof in several steps.
\vskip10pt\par\noindent
{\bf Step 1: Compactness for the cracks.} 
In view of \eqref{cboundn}, using a variant of Helly's theorem (see for instance \cite[Theorem 6.3]{DMT} for the case of Hausdorff converging compact sets), we can find a subsequence (not rabelled) of $(K_n(\cdot))_{n \in \N}$ and an increasing map $t \to K(t)$ such that $K_n(t)$ $\sigma$-converges in $\Omb$ to $K(t)$
for all $t \in [0,T]$.
\vskip10pt\par\noindent
{\bf Step 2: Compactness for the displacements.}
Notice that the sequence $(u_n(t))_{n \in \N}$ is relatively compact in $SBV^p(\Om)$ 
by \eqref{cboundn}. We now want to select a particular limit point of this sequence.
With this aim, let us consider
$$
\vartheta_n(t):= \int_\Om \nabla_\xi f_n(x,\nabla u_n(t))\nabla \dot \psi_n(t)\,dx
\qquad \text{and} \qquad
\vartheta(t):=\limsup_{n \to +\infty} \vartheta_n(t).
$$
Let us see that there exists $u(t) \in SBV^p(\Om)$ such that
\begin{equation}
\label{deftheta}
\vartheta(t)=\int_\Om \nabla_\xi f(x,\nabla u(t))\nabla \dot \psi(t)\,dx
\end{equation}
and
\begin{equation}
\label{weakconvunt}
u_{n_k}(t) \weak u(t)
\qquad
\text{weakly in }SBV^p(\Om)
\end{equation}
for a suitable subsequence $n_k$ depending on $t$. In fact let us consider
a subsequence $n_k$ such that
$$
\vartheta(t)=\lim_{k \to +\infty} \int_\Om \nabla_\xi f(x,\nabla u_{n_k}(t))
\nabla \dot \psi_{n_k}(t)\,dx,
$$
and
$$
u_{n_k}(t) \weak u 
\qquad
\text{weakly in }SBV^p(\Om).
$$
By static equilibrium for $(u_n(t),K_n(t))$ we have that
$$
\int_\Om f_{n_k}(x,\nabla u_{n_k}(t))\,dx \le
\int_\Om f_{n_k}(x,\nabla v(x))\,dx+\int_{H \setminus K_{n_k}(t)} g_n(x,\nu)\,d\hn(x)
$$
for all $v \in AD(\psi_{n_k}(t),H)$ with $K_{n_k}(t) \tsub H$. 
Then by Theorem \ref{stabilitybdrythm} we get that
\begin{equation*}
\label{staticeq}
\int_\Om f(x,\nabla u)\,dx \le
\int_\Om f(x,\nabla v(x))\,dx+\int_{H \setminus K(t)} g(x,\nu)\,d\hn(x)
\end{equation*}
for all $v \in AD(\psi(t),H)$ with $K(t) \tsub H$ and 
\begin{equation*}
\label{convebulk}
\int_\Om f_{n_k}(x,\nabla u_{n_k}(t))\,dx \to
\int_\Om f(x,\nabla u)\,dx.
\end{equation*}
We claim that
\begin{equation}
\label{momenti}
\lim_{k \to +\infty}
\int_\Om \nabla_\xi f_{n_k}(x,\nabla u_{n_k}(t)) \nabla \Phi \,dx
=
\int_\Om  \nabla_\xi f(x,\nabla u)\nabla \Phi\,dx
\end{equation}
for all $\Phi \in W^{1,p}(\Om)$.
This has been done in \cite[Lemma 4.11]{DMFT} in the case of fixed bulk energy, and our proof is just a variant based on the $\Gamma${-}convergence results of Section \ref{gammasec} and on assuption \eqref{fnhyp} which permit to deal with varying energies.
Let us consider $s_j \searrow 0$ and $k_j \to +\infty$: up to a further subsequence for $k_j$ we can assume that
$$
\int_\Om \frac{f(x,\nabla u(x)+s_j \nabla \Phi(x))-f(x,\nabla u(x))}{s_j}\,dx-\frac{1}{j}
\le
\int_\Om \nabla_\xi f_{n_{k_j}}(x, \nabla u_{n_{k_j}}(t)+\tilde s_j \nabla \Phi)
\nabla \Phi\,dx
$$
where $\tilde s_j \in [0,s_j]$. This comes from lower semicontinuity for bulk energies under $\Gamma${-}convergence given by Proposition \ref{lscprop2}, and by Lagrange's Theorem. By Lemma \ref{convmom} we have
$$
\liminf_{j \to +\infty}
\int_\Om \nabla_\xi f_{n_{k_j}}(x, \nabla u_{n_{k_j}}(t)+\tilde s_j \nabla \Phi)
\nabla \Phi\,dx=
\liminf_{j \to +\infty}
\int_\Om \nabla_\xi f_{n_{k_j}}(x, \nabla u_{n_{k_j}}(t))\nabla \Phi\,dx,
$$
so that we get
$$
\int_\Om \nabla_\xi f(x,\nabla u)\nabla \Phi\,dx \le
\liminf_{j \to +\infty}
\int_\Om \nabla_\xi f_{n_{k_j}}(x, \nabla u_{n_{k_j}}(t))\nabla \Phi\,dx.
$$
Changing $\Phi$ with $-\Phi$, we get that 
\eqref{momenti} is proved: setting $u(t):=u$ we deduce that
\eqref{deftheta} and \eqref{weakconvunt} hold.
\vskip10pt\par\noindent
{\bf Step 3: Conclusion.}
Let us consider $t \to (u(t),K(t))$
with $u(t)$ and $K(t)$ defined in Step 2 and Step 1 respectively.
In order to see that $t \to (u(t),K(t))$ is a quasistatic crack evolution we have to check the admissibility condition $u(t) \in AD(\psi(t),K(t))$ for all $t$, and the properties of irreversibility, static equilibrium and energy balance conditions with respect to $f$ and $g$.
\par
As for admissibility, this is guaranteed by \eqref{weakconvunt} and by Proposition \ref{sigmaprop2} which ensures that $\Sg{\psi(t)}{u(t)} \tsub K(t)$.
{\it Irreversibility} is given by construction in Step 1, and {\it static equilibrium} comes from \eqref{staticeq} for $t \in (0,T]$, and by Lemma \ref{gammaconvlem} (where we take $K_n=\emptyset$) for $t=0$. As for {\it energy balance}, we have that static equilibrium implies that (see \cite{DMFT}) for all $t \in [0,T]$
$$
\Es(u(t),K(t)) \ge \Es(u(0),K(0)) +
\int_0^t \int_\Om \nabla_{\xi}f(x,\nabla u(\tau)) \nabla \dot{\psi}(\tau) \,dx\,d\tau.
$$
On the other hand by lower semicontinuity given by Proposition \ref{lscprop2} and by Proposition \ref{lscsigmap} (applied to $g'$ from which $g$ is obtained by restriction) we have for all $t \in [0,T]$
$$
\Es(u(t),K(t)) \le \liminf_{n \to +\infty} \Es_n(u_n(t),K_n(t)),
$$
and by $\Gamma$-convergence given by Lemma \ref{gammaconvlem} (where we take $K_n=\emptyset$)
$$
\Es(u(0),K(0))=\lim_{n \to +\infty}\Es_n(u_n(0),K_n(0)).
$$
Hence we get for all $t \in [0,T]$ (applying also Fatou's Lemma in the limsup version)
\begin{multline}
\label{asymptener}
\nonumber 
\Es(u(t),K(t)) \le \liminf_{n \to +\infty} \Es_n(u_n(t),K_n(t)) \le
\limsup_{n \to +\infty} \Es_n(u_n(t),K_n(t)) \\
=\limsup_{n \to +\infty} \Es_n(u_n(0),K_n(0))+\int_0^t \vartheta_n(s)\,ds
\le
\Es(u(0),K(0))+\int_0^t \vartheta(s)\,ds \\
=\Es(u(0),K(0)) +
\int_0^t \int_\Om \nabla_{\xi}f(x,\nabla u(\tau)) \nabla \dot{\psi}(\tau) \,dx\,d\tau
\le \Es(u(t),K(t)),
\end{multline}
so that we get
$$
\Es(u(t),K(t))=\Es(u(0),K(0))+
\int_0^t \int_\Om \nabla_{\xi}f(x,\nabla u(\tau)) \nabla \dot{\psi}(\tau) \,dx\,d\tau
$$
and
\begin{equation*}
\label{convener}
\lim_{n \to +\infty} \Es_n(u_n(t),K_n(t))=\Es(u(t),K(t)).
\end{equation*}
Finally by lower semicontinuity for the bulk and surface energies under weak convergence for the displacements and $\sigma${-}convergence in $\Omb$ for the cracks, we conclude that
$$
\lim_{n \to +\infty} \Eb_n(u_n(t))=\Eb(u(t))
\qquad \text{and} \qquad
\lim_{n \to +\infty} \Esup_n(K_n(t))=\Esup(K(t)),
$$
so that the theorem is proved.
\end{proof}

\begin{remark}
\label{finalremark}
{\rm
Following the arguments of preceding proof, it turns out that Theorem \ref{compevol}
also holds in the following {\it discretized in time} version, which is closer in spirit to the approach of Francfort and Marigo \cite{FM} to quasistatic crack propagation, and of the subsequent papers on the subject (\cite{AcP}, \cite{Ch}, \cite{DMFT}, \cite{DMT}, \cite{FL}, \cite{GP} and \cite{GP2}). 
\par
Let $0<t^\delta_0<\dots<t^\delta_h=T$ be a subdivision of $[0,T]$ 
with step $\delta>0$, and let $(u^i_{\delta,n},K^i_{\delta,n})$ be such that
$$
(u^i_{\delta,n},K^i_{\delta,n}) \in
\text{ argmin } \{\Eb_n(u)+\Esup_n(K)\,:\, u \in AD(\psi(t^\delta_i),K),\, K^{i-1}_{\delta,n} \tsub K\},
$$
where we set $K^{-1}_{\delta_n}:=\emptyset$.
Let $\delta_n \to 0$, and let $t \to (u_n(t),K_n(t))$ be the {\it discretized in time} evolution defined as
$$
u_n(t):=u^i_{\delta_n,n}, \qquad K_n(t):=K^i_{\delta_n,n}, \qquad
t^i_{\delta_n} \le t < t^i_{\delta_n},
$$
with $u_n(T):=u^h_{\delta_n,n}$ and $K_n(T):=K^h_{\delta_n,n}$.
\par
Then there exists a quasistatic crack growth $t \to (u(t),K(t))$ relative to the energy densities $f$ and $g$ and the boundary datum $\psi$ such that, up to a subsequence (not rabelled), points $(1)$ and $(2)$ of Theorem \ref{compevol} hold.
} 
\end{remark}

\begin{remark}
{\rm 
Notice that for all $t \in [0,T]$ $K_n(t)$ converges to $K(t)$ also in the sense of $\sigma^p$-convergence by Dal Maso, Francfort and Toader \cite{DMFT} (see Section \ref{sigmaconv} just before Corollary \ref{sigmasigmap} for a definition). In fact, by compactness of $\sigma^p$-convergence, up to a further subsequence we have that $K_n(t)$ $\sigma^p$-converges to some $\tilde K(t)$; by Corollary \ref{sigmasigmap}
$\tilde K(t)$ is contained in $K(t)$ so that the pair $(u(t),\tilde K(t))$ is a unilateral minimizer with respect to $f$ and $g$. Following Step $3$ we obtain that $\Esup_n(K_n(t)) \to \Esup(\tilde K(t))$, which together with $\Esup_n(K_n(t)) \to \Esup(K(t))$ implies $K(t) \tilde{=} \tilde K(t)$ for all $t \in [0,T]$.
\par
We conclude that in order to deal with the study of the asymptotic behavior 
of quasistatic crack growths the notion of $\sigma$-convergence and $\sigma^p$-convergence of rectifiable sets are equivalent. 
Notice however that, as pointed out in the Introduction, in order to handle the problem using directly the tool of $\sigma^p$-convergence one would have to prove a Transfer of Jump Sets like our Theorem \ref{transferofjump}, which seems difficult to be derived without any $\Gamma$-convergence argument.
}
\end{remark}

\vskip10pt\par\noindent {\bf Acknowledgments.}
This work began while the authors were visiting the Laboratoire J.-L. Lions of the University of Paris 6 and the L.P.M.T.M. of the University of Paris 13 under the support of the Universit\`a Italo-francese. The authors wish to thank G. Dal Maso and G.A. Francfort for several interesting discussions.

\end{document}